\documentclass[a4paper]{amsart}
\date{\timestamp}
\input epsf.sty
\usepackage[dvips]{graphicx}
\usepackage{verbatim,enumerate}
\usepackage[usenames]{color}

\usepackage{times}
\usepackage{amsthm}
\theoremstyle{plain}
  \newtheorem{theorem}{Theorem}[section]
  \newtheorem*{theorem*}{Theorem}
  \newtheorem{proposition}[theorem]{Proposition}
  \newtheorem{lemma}[theorem]{Lemma}
  \newtheorem{corollary}[theorem]{Corollary}
\theoremstyle{definition}
  \newtheorem{definition}[theorem]{Definition}
\theoremstyle{remark}
  \newtheorem{remark}[theorem]{Remark}
  
\numberwithin{equation}{section}

\newcommand{\vect}[1]{\boldsymbol{#1}}
\newcommand{\R}{\vect{R}}

\newcommand{\op}{\operatorname}

\newcount\minutes
\newcount\scratch

\def\timestamp{%
\scratch=\time
\divide\scratch by 60
\edef\hours{\the\scratch}
\multiply\scratch by 60
\minutes=\time
\advance\minutes by -\scratch
\the\day$.\the \month.\the\year$ at $\,$\hours:\null
\ifnum\minutes< 10 0\fi
\the\minutes}
\title[Inflection Points]{
    Inflection points and double tangents on anti-convex curves
    in the real projective plane
}

\author{Gudlaugur Thorbergsson}
\address[Thorbergsson]{%
  Mathematisches Institut,
Universit\"at zu K\"oln
Weyertal 86 - 90
50931 K\"oln
Germany %
}
\email{gthorbergsson@mi.uni-koeln.de}

\author{Masaaki Umehara}
\address[Umehara]{%
    Department of Mathematics, Graduate School of Science,
    Osaka University,
    Toyonaka, Osaka 560-0043,
    Japan
}
\email{umehara@math.wani.osaka-u.ac.jp}

\medskip

\begin{document}
\begin{abstract}
A simple closed curve $\gamma$ in
the real projective plane $P^2$ is called {\it anti-convex}
if for each point $p$ on the curve, there exists
a line which is transversal to the curve and
meets the curve only at $p$.
We shall prove the relation
$i(\gamma)-2\delta(\gamma)=3$ for anti-convex curves,
where $i(\gamma)$ is the number
of independent (true) inflection points
and $\delta(\gamma)$ the number of
independent double tangents.
This formula is a refinement of the classical M\"obius theorem.
We shall also show that there are three inflection points
on a given anti-convex curve
such  that the tangent lines at
these three inflection points
cross the curve only once.
Our approach is axiomatic and can be applied in other situations.
For example, we prove similar results for curves of
constant width as a corollary.
\end{abstract}

\thanks{
The first author was supported in part by the DFG-Schwerpunkt  
{\it Globale Differentialgeometrie}. \newline
\,\,\,The second author was partly supported by the Grant-in-Aid for 
Scientific Research (B), Japan Society for the Promotion of Science.}
\maketitle

\section*{Introduction}\medskip

Let $P^2$ denote the real projective plane.
We assume curves to be
parameterized and $C^1$-regular.
A simple closed curves in $P^2$ is said to be {\it anti-convex}
or {\it satisfying the Barner condition}
if for each point $p$ on the curve, there exists
a line which is transversal to the curve and
meets the curve only at $p$.
This condition is the  $n=2$ case of a condition
introduced by Barner in \cite{Barner}
for simple closed curves in the real projective space $P^n$  for $n\ge
2$.
An anti-convex curve is automatically  not contractible.

\begin{figure}[htbt]
   \begin{center}
          \includegraphics[width=3.5cm]{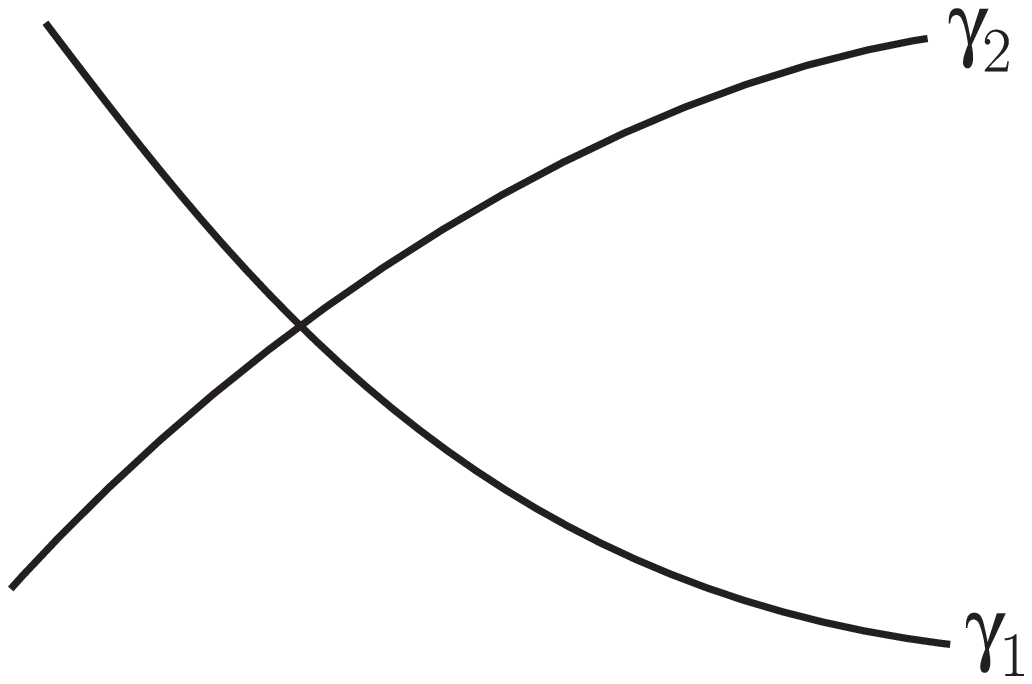}
\qquad
          \includegraphics[width=5cm]{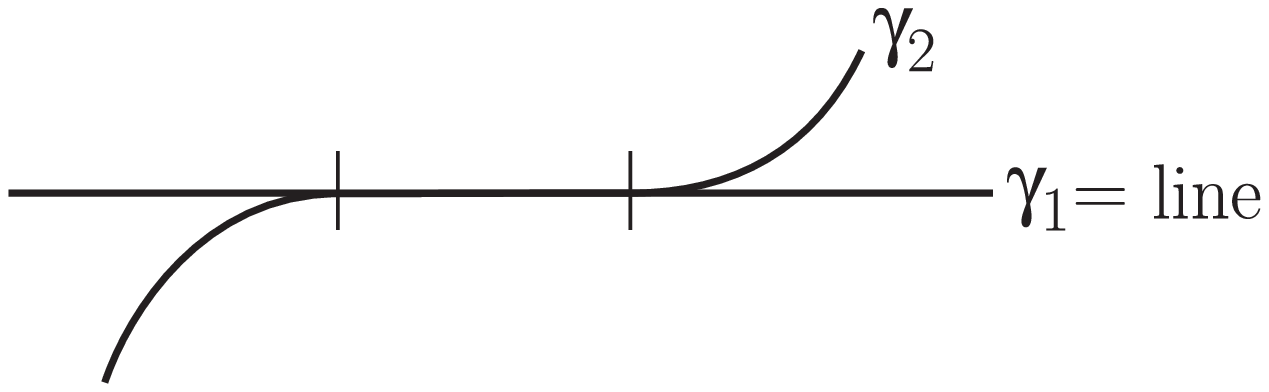}
\caption{Crossing curves}
\label{crossing}
   \end{center}
\end{figure}

Let $\gamma_1$ and $\gamma_2$ be two arcs in
some affine plane $A^2\subset P^2$.
We say that $\gamma_1$ {\it crosses} $\gamma_2$ in
a closed arc $\alpha$
if $\alpha$ is a maximal common  arc of $\gamma_1$ and $\gamma_2$
and there is an open subarc $\tilde\alpha$ of $\gamma_1$
containing $\alpha$ such that the
two
components of $\tilde\alpha-\alpha$ do not lie on the same
side of
$\gamma_2$ (but might not be disjoint from $\gamma_2$).
The arc $\alpha$ can of course
consist of a single point.
If $\gamma_1$ meets $\gamma_2$ transversally in a point $p$, then
$\gamma_1$ of course
crosses $\gamma_2$ in $p$. Examples of crossing curves
are shown in Figure \ref{crossing}.

An inflection point $p$ of a curve $\gamma$ will be called a {\it true
inflection point}
if the tangent line of $\gamma$ at $p$ crosses $\gamma$ in an arc
containing $p$.
Two inflection points are called {\it independent} if
they are not contained
in an arc of $\gamma$ consisting of  true inflection points.
(The inflection points on the curve $\gamma_2$ on the right in Figure 1  
are not indpendent.
On the other hand, the three inflection points in Figure 4 are  
independent.)
We will denote the maximal number of independent true inflection points  
on $\gamma$ by
$i(\gamma)$.

\begin{figure}[htbt]
   \begin{center}
     \includegraphics[width=4cm]{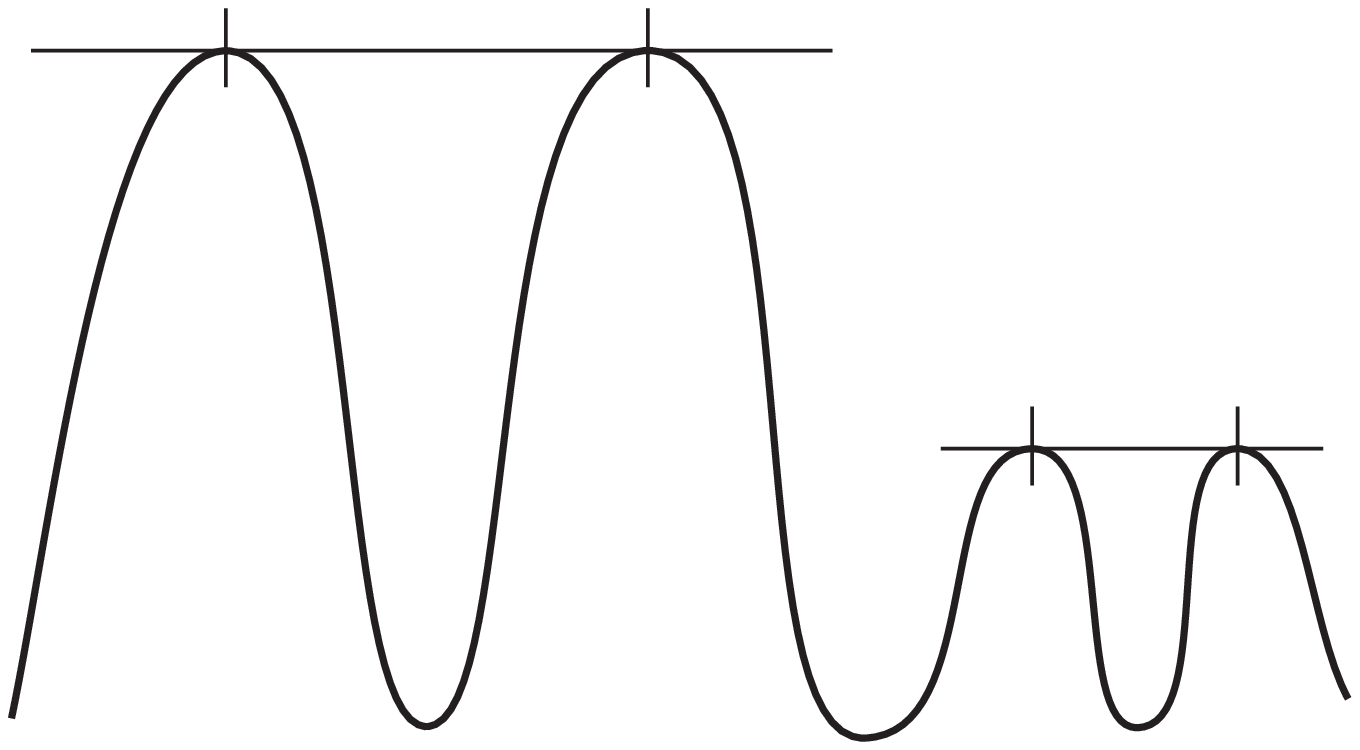}\qquad
\includegraphics[width=4cm]{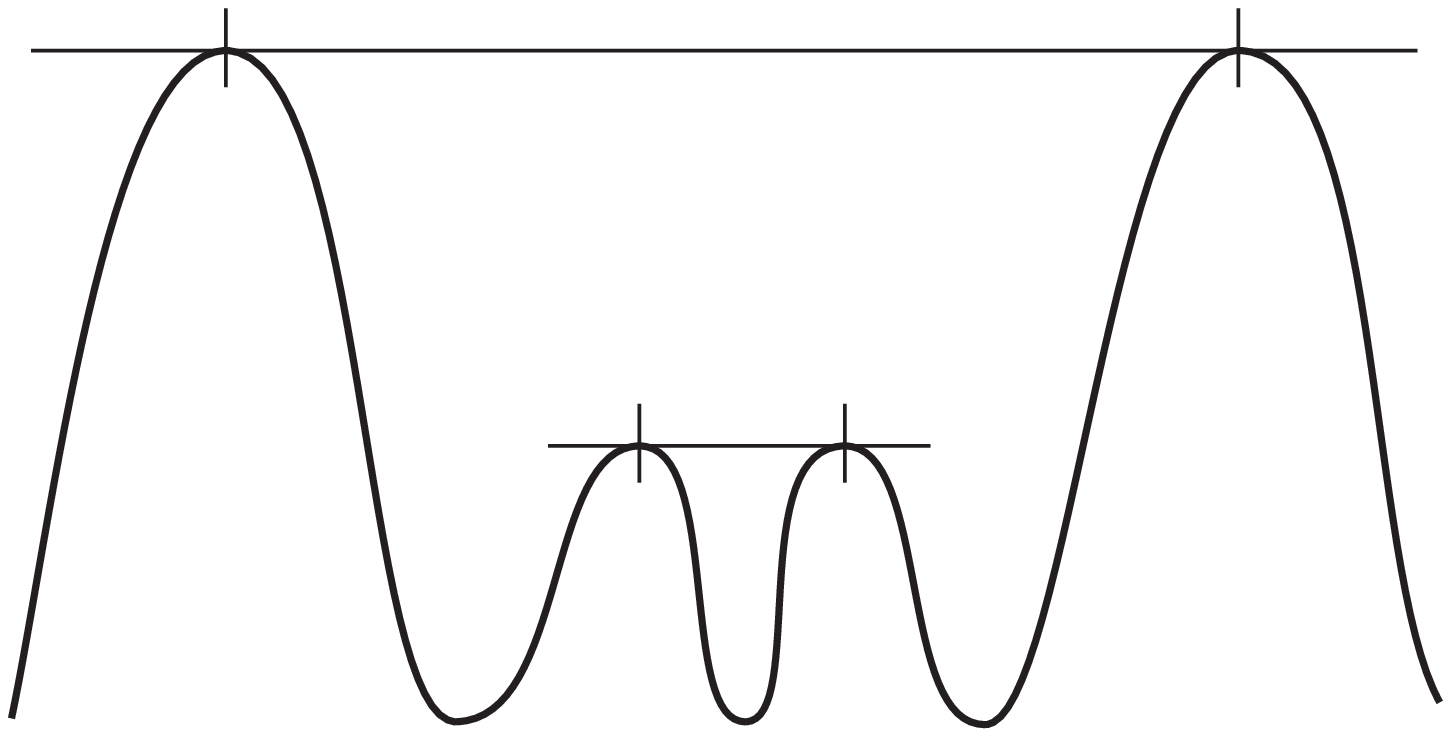}
\caption{Two types of independent double tangent}
\label{independent}
   \end{center}
\end{figure}

\begin{figure}[htbt]
   \begin{center}
     \includegraphics[width=5cm]{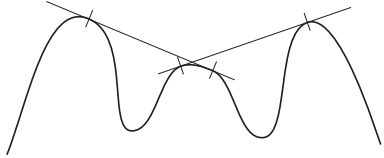}
\caption{Dependent double tangents}
\label{non-indep}
   \end{center}
\end{figure}

  A {\it double tangent} of a curve $\gamma$ is roughly speaking a
line $L$ that is tangent to $\gamma$  at the endpoints of a nontrivial  
arc $\alpha$
of $\gamma$ contained in an affine plane $A^2\subset P^2$ in such a way  
that $\alpha$
is locally around its endpoints on the same side of $L\cap A^2$.
(A precise definition will be given in Section 4.)
We call  $\alpha$ a {\it double tangent
arc}.  A set of double
tangent arcs $\alpha_1,\dots\alpha_k$ is said to be {\it independent}
if  any two of the arcs
   are either disjoint or one is a subarc of the other; see
Figures~\ref{independent} and \ref{non-indep}.

We will denote the number of elements in a maximal  set of
independent double tangent arcs
by $\delta(\gamma)$. It will follow from Theorem A, which we now state,
that $\delta(\gamma)$ is
independent of the choice of a maximal set of independent double  
tangent arcs on
$\gamma$.

\medskip
\noindent{\bf Theorem A.} \, {\it Let $\gamma$
be a $C^1$-regular
anti-convex curve in $P^2$ which is not a line.
If the number $i(\gamma)$ of independent
true inflection points on $\gamma$ is finite, then so is
the  number $\delta(\gamma)$ of elements in a maximal set of
independent double  tangents,
and
$$
i(\gamma)-2\delta(\gamma)=3 \leqno (*)
$$
holds. In particular, the number $\delta(\gamma)$
does not depend on the choice of  a maximal set of
independent  double tangents if $i(\gamma)$ is finite.}\medskip
\medskip

Formula $(*)$ is reminiscent of the Bose formula  for simple
closed curves in the Euclidean plane  saying that $s-t=2$,  where
$s$ is the number of inscribed osculating circles
and $t$ is the number of triple tangent inscribed circles. This formula
was proved for convex curves by Bose in \cite{Bose}
and in the general case by
Haupt in \cite{Haupt1}.
Our method to prove Theorem A will be similar to the one used by the
second author to prove the Bose formula in \cite{Umehara}.
The authors do not know whether  formula $(*)$ holds
for non-contractible simple closed curves which are not necessarily   
anti-convex.

There is a well-known formula
for generic closed curves in the affine plane $A^2$  due to
Fabricius-Bjerre
relating the numbers of double points, inflection points, and double  
tangents; see
\cite{FB}.
When the curves have no inflection points, Ozawa \cite{Ozawa} gave
a sharp upper bound on the number of double tangents. Formulas for
real algebraic curves in $P^2$ go at least back to Klein; see the paper
\cite{Wa} of Wall.\medskip

We will also prove the following theorem.\medskip

\noindent{\bf Theorem B.} \,  {\it Let $\gamma$
be a $C^2$-regular anti-convex curve in $P^2$ which is not a line.
Then $\gamma$ has at least three inflection points
with the property that the tangent lines at these inflection points
cross $\gamma$  only once.}
\medskip

\begin{figure}[htbt]
   \begin{center}
          \includegraphics[width=7cm]{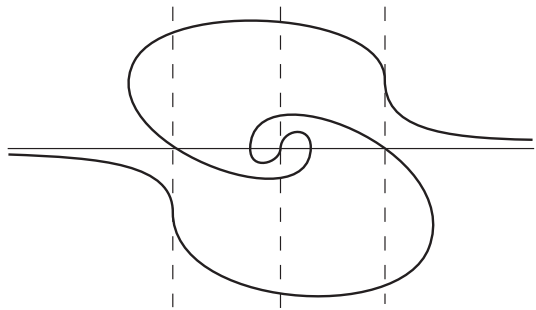}
\caption{A simple closed curve with no clean inflection points}
\label{spiral}
   \end{center}
\end{figure}

The theorem is optimal.
An inflection point $p$ is called {\it clean} if the
tangent line at $p$ meets the curve in a connected set.
A clean inflection point is a typical example of an
inflection point as in Theorem B.
The noncontractible branch of a
regular cubic in $P^2$ has three clean inflection points.
   M\"obius proved that a simple closed
noncontractable curve in $P^2$ has at least three (true) inflection
points.
Several proofs this result are known; see \cite{A.K},  \cite{Haupt2},
and \cite{Sasaki}.
One can show with examples that none of these has to be
a clean inflection point; see Figure \ref{spiral}.

A similar result is proved in \cite{TU2} and \cite{TU3}
for clean sextactic points on a strictly convex curve in the affine  
plane.
It says that such a curve has
{\it  three inscribed osculating conics and
three circumscribed osculating conics}.
It should also be remarked that the Tennis Ball
Theorem (\cite{A1} and \cite{A2}), the
theorem of Segre on space curves in \cite{Segre},
and the refinement of the  Four-Vertex Theorem in
\cite{TU} can be considered as
generalizations of the M\"obius Theorem; see \cite{TU}.

In the  proof of Theorem B
we use an approach that goes back to H. Kneser's proof of
the four vertex theorem; see \cite{K.H}, \cite{Umehara},
and also  \cite{TU}, \cite{Th}.
(A further development of this approach is also crucial in the proof of  
Theorem A.)

\begin{figure}[h]
   \begin{center}
          \includegraphics[width=5cm]{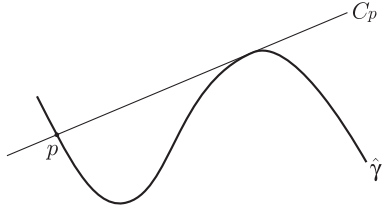}
\caption{The limiting great circle}
\label{limiting}
   \end{center}
\end{figure}

\medskip
The theorems will be proved in later sections.
Here we would like to explain some of the basic ideas in the proofs.
   Let $\hat \pi:S^2\to P^2$ be the universal covering of $P^2$. Since
$\gamma$ is not contractible, it lifts to a simple closed curve
$\hat\gamma$ that double covers $\gamma$. There is through every
point $p$ on $\hat\gamma$ a great circle $\hat L_p$ on $S^2$
(which is the double cover of the line $L_p$) that only
meets $\hat\gamma$ in $p$ and the antipodal point $T(p)=-p$.
The parametrization of $\hat\gamma$ and the orientation of $S^2$
give us a tangent and normal vector field along $\hat\gamma$. We
will assume that
the normal direction points to the left side of the curve. We define a
positive rotation direction along
the curve by rotating the normal vector towards the tangent vector.
Notice that the  positive rotation direction
is the clockwise direction.
Let us now  rotate the circle $\hat L_p$ around $p$ as far as
possible in the positive direction through circles which only meet
$\hat\gamma$ in $p$ and $T(p)$.
We denote
the limiting great circle by $C_p$.
There are two possibilities. The first is that $C_p$ only meets
$\hat\gamma$
in one component. Then $p$ is a clean inflection point.
The other possibility is that  $C_p$ meets  $\hat\gamma$
in more than one component; see Figure \ref{limiting}.
   In this case  $p$ may or may not be an inflection point,
but it is of course not a clean  inflection point.
We define a closed subset $F(p)$ by setting
\begin{equation}\label{closed}
F(p)=C_p\cap \hat \gamma.
\end{equation}

\medskip
We identify $S^1$ with the image of the curve $\hat\gamma$
and introduce on $S^1$ a cyclic order that agrees with the
orientation
of the curve. We will first assume that no line meets $\gamma$ in
infinitely
many points and then discuss the general case.
If $p$ in $S^1$ is not an inflection point, we let
$\delta$ denote the distance from $p$ to the next point $q\in F(p)$ in
$(p,Tp)$,
where $(a,b)$ denotes the interval from $a\in S^1$ to $b\in S^1$ with
respect
to the cyclic order of $S^1$ and $F(p)$ is defined in equation
\ref{closed}.
Let $p_1$ be the midpoint of the interval $[p,q]$.  The subset
$F({p_1})$
lies in the interval $[p,q]\cup[Tp,Tq]$.
If $p_1$ is not a clean inflection point
we let $\delta_1$ denote the distance to the point $q_1$
closest to $p_1$ in  $F({p_1})\cap (p_1,Tp_1)$.
Notice that $\delta_1\le\delta/2$. Iterating this process, we either
arrive at a point $p_n$
which is a
clean inflection point, or we get a sequence $(p_n)$ that converges to
a clean
inflection point. As we will see in  Section 2, this approach
leads to the existence of
at least three inflection points.
In the proof of Theorem B we only use a few axiomatic properties of
the family  $\{F(p)\}_{p\in S^1}$ of closed subsets in $S^1$. It
can therefore be applied to
different situations.

In Section 5,
we apply the method to convex curves of constant width.

\begin{figure}[h]
   \begin{center}
          \includegraphics[width=5.5cm]{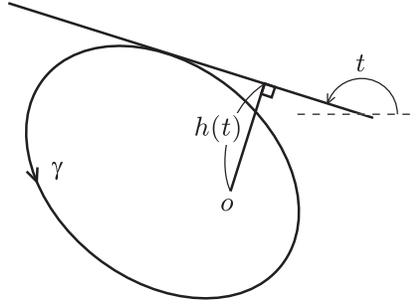}
\caption{The supporting function}
\label{support}
   \end{center}
\end{figure}

Let $\gamma$ be a strictly convex curve in $\R^2$.
For each $t\in [0,2\pi)$,
there is a unique tangent line $L(t)$ of the
curve which makes angle $t$ with the $x$-axis. Let $h(t)$ be the
distance
between a fixed point
$o$  in the open domain bounded by $\gamma$ and the line $L(t)$; see
Figure \ref{support}.  Note that  $t$ gives a
parametrization of the strictly convex curve $\gamma$
which we will use from now on.
   The function $h$ is called the {\it supporting
function of
the curve $\gamma$ with respect to $o$.}
A strictly convex curve has {\it constant width} $d$ if and only if
$h(t)+h(t+\pi)=d$ holds.

We now  fix a curve $\gamma$ of constant width $d$.
For each point $p$ on the curve, there exists a
unique circle $\Gamma_p$
of width $d$ such that  $\Gamma_p$ is tangent to $\gamma$ at $p$,
that is $\Gamma_p$ and $\gamma$ meet at $p$ with multiplicity two.
Since $\Gamma_p$ is the best approximation of $\gamma$ at $p$ among the
circles of width $d$, we call $\Gamma_p$ the {\it osculating
$d$-circle at $p$.}
Generically, the osculating $d$-circle of $\gamma$ at $p$ does not  
cross $\gamma$ at $p$.

We will  prove the following theorem in Section 5.

\medskip
\noindent{\bf Theorem C.} \, {\it Let $\gamma$ be
a $C^3$-regular strictly convex curve of constant width $d$.
Then there exist at least three osculating $d$-circles
which cross  $\gamma$ exactly twice, both times
tangentially.
Moreover, these three circles coincide with
the osculating circles (in the usual sense)
at each of their crossing points on $\gamma$.
}
\medskip

\begin{figure}[htb]
\begin{center}
\includegraphics[height=6cm]{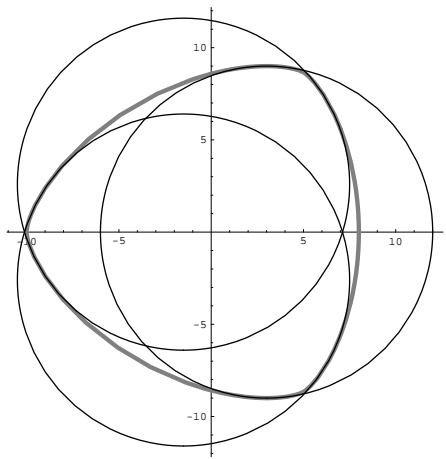}
\caption{The three osculating circles}
\label{width}
\end{center}
\end{figure}

The above theorem is a refinement of the fact that
there are six distinct points
on $\gamma$ whose osculating circles have
radius $d/2$.
  (Basic properties of curves of constant width can be found
in \cite{Y}.)
In Figure \ref{width} we indicate the three osculating circles
of diameter $d$ of the curve of constant width whose
supporting function is $(d/2)+\sin 3t$.

We will also prove a  formula analogous to the one in Theorem A for
curves
of constant width in Section 5.
\bigskip

\section{Intrinsic line systems} \medskip

In this section, we shall derive some basic properties
of the family of closed subsets $\{F(p)\}_{p\in S^1}$
defined in equation \eqref{closed} in the introduction.
We shall then
use these properties to define
what we  will call an
\lq intrinsic line system\rq.

Let $\gamma:P^1\to P^2$ be a $C^1$-regular
anti-convex curve in $P^2$,
where $P^1$ is a closed circle considered as a projective line.
We assume that the image of $\gamma$ is not a line in $P^2$.
Let
$\hat \pi:S^2\to P^2$ and $\pi:S^1\to P^1$
be the canonical covering projections.
Then there exists a simple closed curve
$\hat \gamma:S^1\to S^2$
such that
$$
\hat \pi\circ \hat \gamma=\gamma\circ \pi.
$$
Moreover, for each point $p$ on $\hat \gamma$, there exists
a great circle $\hat L_p$ on $S^2$ such that
$\hat \pi(\hat L_p)=L_{\pi(p)}$.
By rotating $\hat L_p$  in the clockwise direction through great
circles
that only meet $\hat\gamma$ in $p$ and the antipodal point $Tp$, we
arrive at the limiting great
circle $C_p$ as in the introduction.
Let $D_{\hat \gamma}$ be the domain
on the left hand side of $\hat \gamma$.
We  orient $\hat L_p$ such that
it passes into $D_{\hat \gamma}$ after going through $p$. The
orientation of
the great circle $\hat L_p$ induces an orientation on the limiting
great circle $C_p$.

If $C$ is an oriented great circle, we denote by $H^+(C)$ (resp.
$H^-(C)$)
the closed hemisphere on
the left (resp. right) hand side of $C$.

\begin{figure}[htb]
   \begin{center}
         \includegraphics[width=4.5cm]{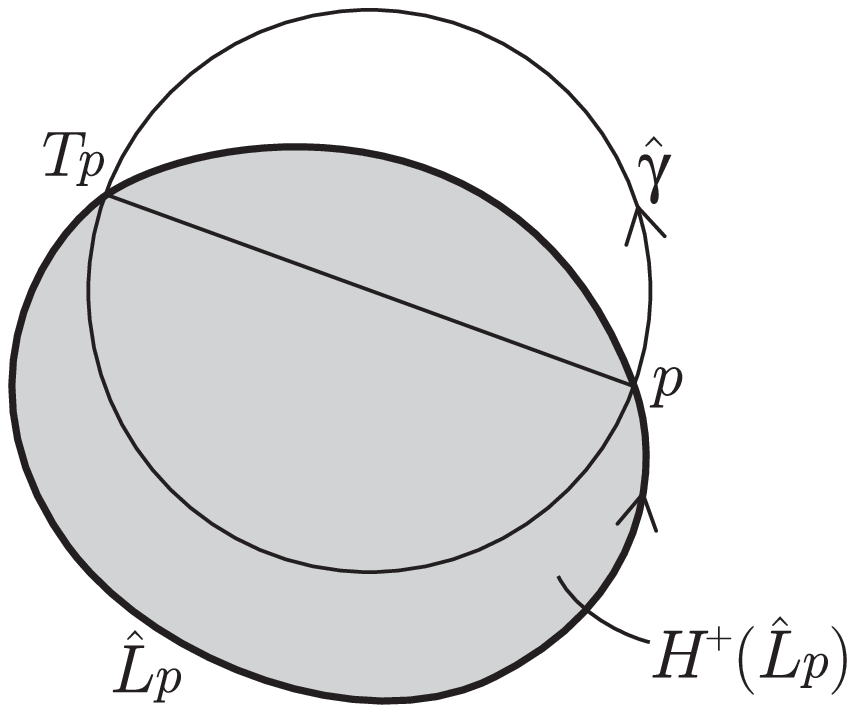}
\caption{}
\label{pic1}
   \end{center}
\end{figure}

By applying a suitable diffeomorphism to $S^2$, we can
map $\hat \gamma$ onto the equator and $D_{\hat\gamma}$
on  the upper hemisphere.
If we compose this with  the stereographic projection into
   the plane, $\hat\gamma$ and
$H^+(\hat L_p)$ look as in Figure \ref{pic1}.

Though $\hat \gamma$ may not be star-shaped in general,
we shall frequently use this  kind of sketches of $\hat \gamma$
to simplify the figures.

The following assertion is obvious.

\begin{proposition}\label{1.1}
{\it The arc of $\hat \gamma:S^1\to S^2$ from $p$ to $Tp$
$($resp.~from $Tp$ to $p)$ lies in $H^-(\hat L_p)$
$($resp.~$H^+(\hat L_p))$. }
\end{proposition}

\begin{proposition}\label{1.2}
The limiting great circle $C_p$ has the following properties.
\begin{enumerate}
\item[(a)]
The arc of $\hat \gamma$ from $p$ to $Tp$
$($resp.~from $Tp$ to $p)$ lies in $H^-(C_p)$
$($resp.~$H^+(C_p))$.
\item[(b)]
The set $F(p)$ has at least three connected components,
if $C_p$ is not the tangent line of $\hat\gamma$ at $p$.
\end{enumerate}
\end{proposition}

\begin{figure}[htb]
   \begin{center}
         \includegraphics[width=3.7cm]{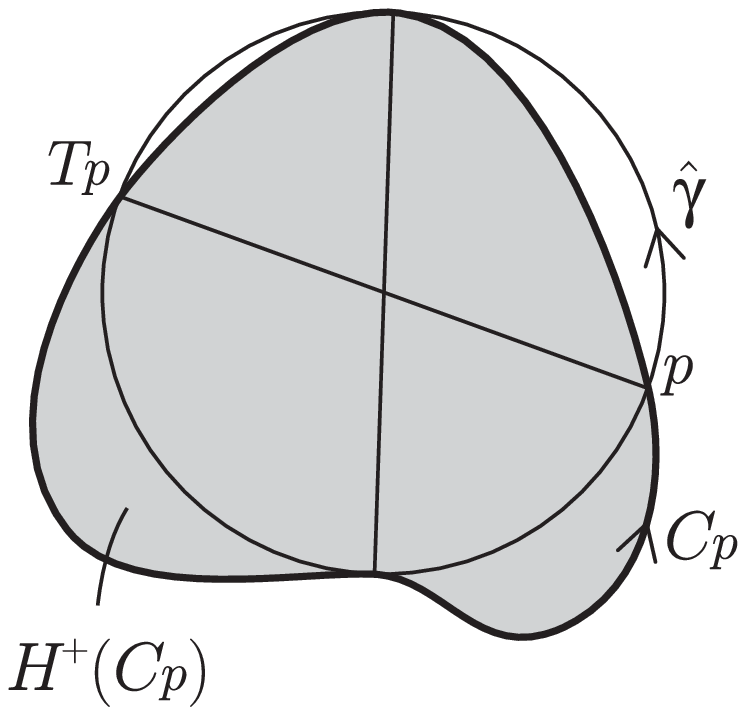}
\caption{}
\label{pic2}
   \end{center}
\end{figure}

\begin{proof}
Since $C_p$ is the limit of circles like $\hat L_p$,
the property in (a) follows from Proposition \ref{1.1}.
To prove (b), we
suppose that  $C_p$ is not a tangent line of $\hat\gamma$ at $p\in S^1$.
Then $C_p$ meets $\gamma$ transversally at $p$ and $Tp$.
Hence if $C_p$ only meets $\hat \gamma$ in these two points,
one can rotate it slightly in positive direction through curves that
are transversal to $\hat\gamma$ in $p$ and $Tp$ and only meet
$\hat\gamma$
in these two points. This contradicts the definition of $C_p$.
Thus there exits a point $q$ in $F(p)=C_p\cap \hat\gamma$ which
is distinct from both $p$ and $Tp$.
Since $\hat \gamma$ is not a great circle, $p$ and $Tp$
belong to different connected components of $F(p)$.
Since both $C_p$ and $\hat \gamma$ are symmetric with respect to
$T$,
it follows that $C_p$ is neither a tangent line at $p$ nor at $Tp$.
If $q$ is in the same connected component of $F(p)$ as $p$
(or $Tp$),
$C_p$ contains the  segment of $\hat\gamma$ between $p$ and
$q$ (or $Tp$ and $q$), which implies that $C_p$ must be the
tangent line at $p$ (resp. $Tp$), a contradiction.
\end{proof}

Conversely, we have the following

\begin{proposition}\label{1.3}
If a great circle $C$ through $p$ and $Tp$ satisfies the
following two properties,
then $C$ coincides with $C_p$.
\begin{enumerate}
\item[(a)]
The arc of $\hat \gamma$ from $p$ to $Tp$
$($resp.~from $Tp$ to $p)$ lies in $H^-(C)$
$($resp.~$H^+(C))$.
\item[(b)]
$C$ is tangent to $\hat\gamma$ at all points in $C\cap\hat \gamma$
different from $p$ and $Tp$ and if $C$ is not tangent to $\hat\gamma$
at $p$ and $Tp$, then
$C\cap\hat \gamma$ contains a point different from $p$ and $Tp$.
\end{enumerate}
\end{proposition}

\begin{proof}
Since $C$ is tangent to $\hat\gamma$ at all points in $C\cap\hat
\gamma$
different from $p$ and $Tp$,
we can rotate $C$ slightly in negative direction into a great circle
which meets $\hat\gamma$ transversally in $p$ and $Tp$ and does not
have any
further points with it in common. It now follows from the definition of
$C_p$ that
$C=C_p$.
\end{proof}

We will denote by $F_0(p)$  the connected component of $F(p)=C_p\cap
\hat \gamma$
containing $p$ for each point $p$ on $S^1$.

\begin{proposition}\label{1.4}
Suppose that  $\gamma:P^1\to P^2$ is an anti-convex curve which is
not a line and
meets a line in $P^2$ in at most
finitely many connected components.
Then the corresponding family  $\{F(p)\}_{p\in S^1}$ of
subsets of $S^1$
satisfies the following properties:
\item{(L1)} $p\in F(p)$.
\item{(L2)} $F(p)$ is a closed proper subset of $S^1$
and has finitely many connected components.
\item{(L3)} If $q\in F(p)$, then $Tq\in F(p)$ where
$T:S^1\to S^1$ is the restriction of
the antipodal map on $S^2$ to $\hat\gamma$.
\item{(L4)} Suppose $p'\in F(p)$ and $q'\in F(q)$ satisfy
$$
p\le q\le p'\le q'(\le Tp)
$$
or
$$
p\ge q\ge p'\ge q'(\ge Tp),
$$
where $\ge$ and $\le$ are the cyclic order of $S^1$.
Then $F(p)=F(q)$.
\item{(L5)}  If $\pi(F(p))=\pi(F_0(p))$, then
$\pi(F(Tp))\ne \pi(F_0(Tp))$ where
$\pi:S^1\to P^1$ denotes the canonical projection.
\item{(L6)} $q\in F_0(p)$ if and
only if $F(p)=F(q)$.
\item{(L7)}  Let $(p_k)$ be a sequence in $S^1$
that converges to an element
   $p$ in $S^1$, and  let $(s_k)$  be another sequence in $S^1$ such that
$s_k\in
F(p_k)$ and $\lim s_k=s$.
Then $s\in F(p)$.
\end{proposition}

\begin{proof}
(L1) is obvious. (L2) is a trivial consequence of
the assumption that $\gamma$ and a line meet in
at most finitely many connected components.
(L3) follows from the fact that $\hat\gamma$ and $\hat L_p$
are both symmetric with respect to the antipodal map $T$.

We  now prove (L4).
If $C_p$ and $C_q$ are  great circles
which meet in two points which are not antipodal,
then $C_p$ must be equal to $C_q$.
Suppose $p'\in F(p)$ and $q'\in F(q)$ and
$
p\le q\le p'\le q'(\le Tp)
$
or
$
p\ge q\ge p'\ge q'(\ge Tp)
$
holds. Then the subarc of $C_q$ between $q$ and $q'$ must meet
$C_p$ twice.  One is between $p$ and $p'$, and the other is
between $p'$ and $Tp$ on $C_p$.
(See Figure \ref{two} for the case $p\le q\le p'\le q'$.)
Thus $C_p=C_q$ holds.

\begin{figure}[htb]
   \begin{center}
          \includegraphics[width=5.3cm]{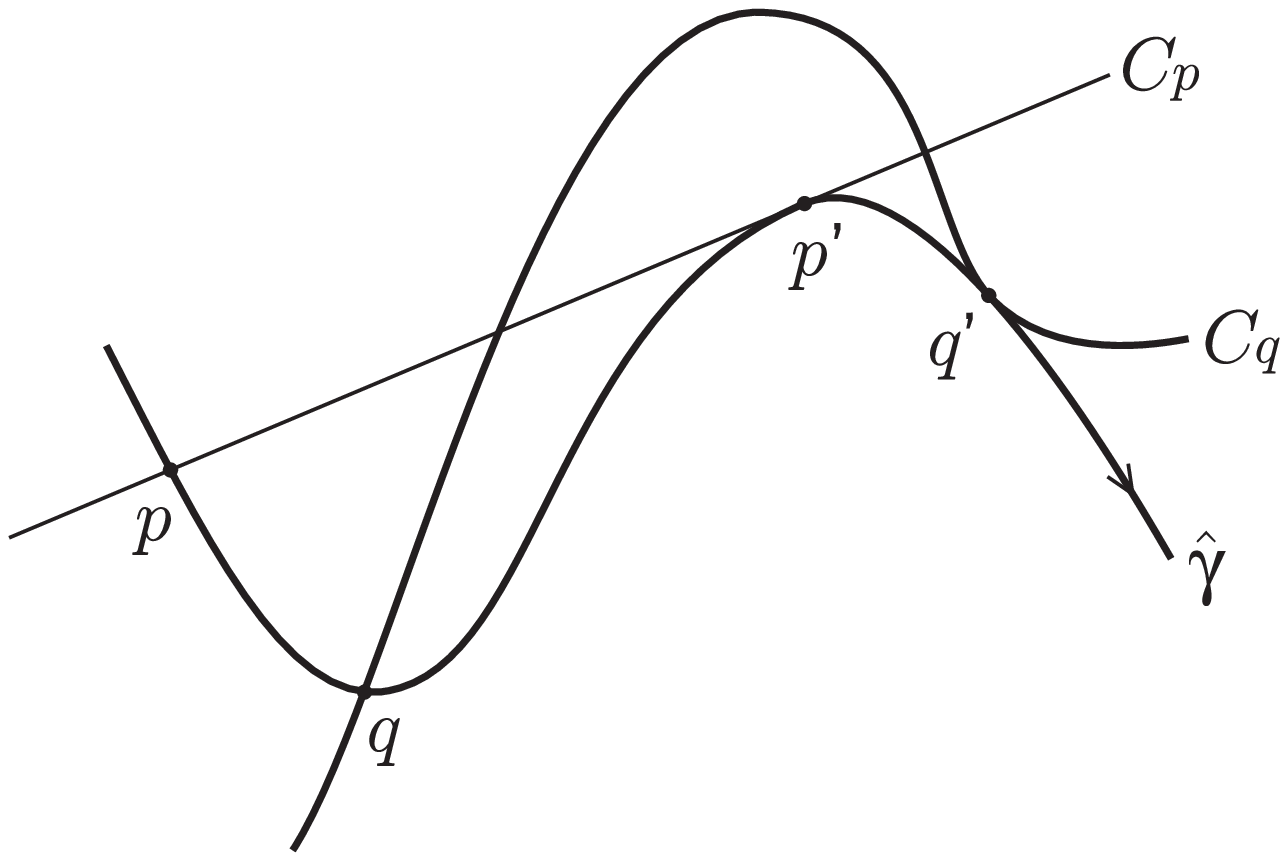}
\caption{}
\label{two}
   \end{center}
\end{figure}

Now we  prove (L5). If $\pi(F(p))=\pi(F_0(p))$, then
$F(p)$ consists of two connected components. By Proposition \ref{1.2}
(b),
$C_p$ is a tangent line at $p$.  The great circle $C_{Tp}$ coincides
with the great circle which we get by rotating $C_p$ in negative
direction
through great circles meeting $\hat\gamma$ only in $p$ and $Tp$ until
it hits
$\hat\gamma$. The great circles $C_p$ and $C_{Tp}$ cannot coincide since
$\gamma$ is not a line.  It follows
that $C_{Tp}$ is not tangent to $\hat\gamma$
at $p$ and hence also not at $Tp$. By Proposition \ref{1.2} (b),
$F(Tp)$ contains
at least three components, two of which consist of $p$ and $Tp$ since
the
intersection between $C_{Tp}$ and $\hat\gamma$ is transversal in these
points.
Hence $\pi(F(Tp))$ is not connected and we see that  $\pi(F(Tp))\ne
\pi(F_0(Tp))$.

We now prove (L6).
Suppose $q\in F_0(p)$. We may assume that
$q\ne p$. Then $F_0(p)$ is a closed interval and
   $C_p$ must be the tangent line both at $p$ and
   $q$. It follows that $C_p$
must be equal to  the great circle $C_q$ by
Proposition \ref{1.3}. This implies $F(p)=F(q)$.
Now we assume that $F(p)=F(q)$. We let $A$ denote
the set of points $r$ in $F(p)=F(q)$ such that the tangent great circle
of $\hat \gamma$ in $r$ contains $F(p)=F(q)$ and $r$ is not a true  
inflection
point.
Let $B$ denote the complement of
$A$ in $F(p)=F(q)$. By Proposition \ref{1.2}  the set $B$
coincides with $F_0(p)\cup T(F_0(p))=F_0(q)\cup T(F_0(q))$.
   Now note that a set $T(F_0(r))$
cannot coincides with a set $F_0(s)$ for any $r$ and $s$ in $S^1$
since the curve $\hat\gamma$ crosses  $C_r$ from right to left in
$F_0(r)$
and  $C_s$ from left to right in $F_0(s)$; see Figure \ref{two-infl}.

Finally we prove (L7).
We may assume that $s$ is neither $p$ nor $Tp$.
After replacing $(p_k)$ by a subsequence if necessary, we
may also assume that $C_{p_k}$ converges to
a great circle $C$.
Since $C_{p_k}$ satisfies properties (a) and (b) in
Proposition \ref{1.3} for all $k$, so does $C$, and it
follows that   $C=C_p$ holds.  Hence $s\in F(p)$.
\end{proof}

\begin{remark}
We will call a family $\{F(p)\}_{p\in S^1}$ of closed subsets of $S^1$
   an {\it intrinsic line system}
if it satisfies properties (L1) -- (L7) in Proposition~\ref{1.4}.
This is an analogue of the somewhat simpler intrinsic circle
systems,
see \cite{Umehara} and \cite{TU}, which are e.g.~useful in
proving the existence of two inscribed (resp. circumscribed)  osculating
circles
of a given simple closed $C^2$-regular curve in the Euclidean
plane.
\end{remark}

\begin{figure}[htb]
   \begin{center}
         \includegraphics[width=7cm]{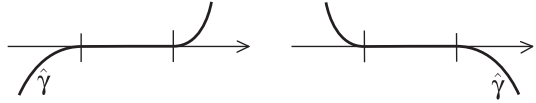}
\caption{Negatve inflection and positive inflection.}
\label{two-infl}
   \end{center}
\end{figure}

An infection point of a curve $\hat\gamma$ in
is called {\it positive} if
the tangent great circle crosses $\hat\gamma$ from right to
left, and {\it negative} if the tangent great circle
crosses $\hat\gamma$ from left to
right.

\begin{definition}\label{1.5} Let $\{F(p)\}_{p\in S^1}$ be an
intrinsic line system.
A point $p\in S^1$ satisfying
$$
\pi(F(p))=\pi(F_0(p))
\qquad ({\rm resp.}\,\, \pi(F(Tp))=\pi(F_0(Tp))),
$$
   is  called a {\it positive clean
infection point} (resp.  a {\it negative clean inflection point}).
\end{definition}

A positive (resp.~negative) clean infection point of $\hat\gamma$ is a  
positive
(negative) inflection point by definition.
Since the sign of a clean inflection point is reversed by the antipodal  
map,
the notion is meaningful for $\hat \gamma$ but not for $\gamma$.
\bigskip

\section{Clean inflection points}\medskip

In this section we prove Theorem B in
the introduction.
The crucial point is that we only use
properties (L1)--(L6) of intrinsic line systems to prove
the theorem under the assumption that
$\gamma$  meets a line in at most finitely many components. It
is only in the last step where we remove this assumption that we
use special properties of curves in $P^2$.

\begin{lemma}\label{2.1} Let $p\in S^1$.
Suppose that $q\in F(p)\cap (p,Tp)$.
Let $r$ be a point in $(p,q)$.
Suppose that $r$ is not contained in  $F_0(p)$.
  Then $$\pi(F(r))\subset\pi((p,q))$$ holds.
\end{lemma}

\begin{proof}
Suppose $\pi(F(r))$ contains an element $a\not\in \pi((p,q))$. Let
$\{\hat a^+,\hat a^-\}$  be the preimage of $a$ under $\pi$.
Without loss of generality, we may assume that $\hat a^+\in (p,Tp]$.
Since $a\not\in \pi((p,q))$, we have $\hat a^+\in (q,Tp]$.
Hence we have the inequality
$$
p < r< q \le {\hat a}^+ \le Tp.
$$
By (L4), we have $F(p)=F(r)$.
In particular $r\in F_0(p)$ by (L6), which is a contradiction.
\end{proof}

With similar arguments we can prove the following lemma.

\begin{lemma}\label{2.2}
Let $p\in S^1$.
Suppose that $q\in F(p)\cap (Tp,p)$.
Let $r$ be a point in $(q,p)$.
Suppose that $r$ is not contained in  $F_0(p)$.
  Then $$\pi(F(r))\subset\pi((q,p))$$ holds.
\end{lemma}

Next we prove the following lemma.

\begin{lemma}\label{2.3}
  Let $p\in S^1$.
Suppose that $q\in F(p)\cap (p,Tp)$
  and $(p,q)\cap F_0(p)=\emptyset$.
Let  $r$ be the midpoint of $(p,q)$.
Then at least one of the following three cases occurs:
\item{(i)} $r$ is a positive clean inflection point.
\item{(ii)} There exist $p_1,q_1\in F(r)\cap (r,q)$
such that $p_1\in F_0({r})$ and $(p_1,q_1)\cap F_0(r)=\emptyset$.
\item{(iii)} There exist $p_1,q_1\in F(r)\cap (p,r)$ such
that $p_1\in F_0({r})$ and $(q_1,p_1)\cap F_0(r)=\emptyset$.
\end{lemma}

\begin{proof}
Assume that $r$ is not a positive clean inflection point.
Then there exists a point $b\in \pi(F(r))$, such that
$b\not\in \pi(F_0(r))$.
Let $\{q_1,Tq_1\}$ be the points such that $\pi(q_1)=b$.
Since $(p,q)\cap F_0(p)=\emptyset$, we have $r\not\in F_0(p)$.
Thus by Lemma \ref{2.1}, we have $b\in \pi(F(r))\subset\pi((p,q))$.
So we may assume that $q_1\in (p,q)$ without loss of generality.
Since $b\not\in \pi(F_0(r))$, we have $q_1\not\in F_0(r)$.
There are two possibilities, one being $q_1\in (r,q)$ and the
other being $q_1\in (p,r)$.

First, we consider the case $q_1\in (r,q)$.
Since $F_0(r)$ is a proper subset of $S^1$,
it is a linearly ordered set with respect to
the restriction of the cyclic order of $S^1$
and one can define its supremum and infimum.
We set
$$
p_1:={\rm sup}(F_0(r)).
$$
Since $F_0(r)\subset(p,q)$ and $r\in F_0(r)$,
it holds that $p_1\in [r,q]$.
On the other hand, since $q_1\not\in F_0(r)$ and $q_1\in (r,q)$,
we have
$$
r\le p_1<q_1<q.
$$
This is case (ii).

Next, we consider the case $q_1\in (p,r)$.
We set
$$
p_1:={\rm inf}(F_0(r)).
$$
Since $F_0(r)\subset(p,q)$ and $r\in F_0(r)$,
it holds that $p_1\in [p,r]$.
On the other hand, since $q_1\not\in F_0(r)$ and $q_1\in (p,r)$,
we have
$$
r\ge p_1>q_1>p.
$$
This is case (iii).
\end{proof}

Similarly we get the following lemma.

\begin{lemma} \label{2.4} Let $p\in S^1$.
Suppose that $q\in F(p)\cap (Tp,p)$
  and $(q,p)\cap F_0(p)=\emptyset$.
Let  $r$ be the midpoint of $(q,p)$.
Then at least one of the following three cases occurs:
\item{(i)} $r$ is a positive clean inflection point.
\item{(ii)} There exist $p_1,q_1\in F(r)\cap (r,p)$
such that $p_1\in F_0({r})$ and $(p_1,q_1)\cap F_0(r)=\emptyset$.
\item{(iii)} There exist $p_1,q_1\in F(r)\cap (q,r)$ such
that $p_1\in F_0({r})$ and $(q_1,p_1)\cap F_0(r)=\emptyset$.
\end{lemma}

We will use Lemma \ref{2.3} and Lemma \ref{2.4} to prove the following  
proposition.

\begin{proposition}\label{2.5}
Let $p\in S^1$. Suppose that $q\in F(p)\cap (p,Tp)$
and $(p,q)\cap F_0(p)=\emptyset$.
Then there exists a positive clean inflection point $s$ in $(p,q)$
such that $\pi(F(s))\subset \pi((p,q))$.
\end{proposition}

\begin{proof}
Suppose that there are no positive clean inflection points
in $(p,q)$. Let $\delta$ be the length of the interval $(p,q)$.
Let $r$ denote the midpoint of the interval $(p,q)$.
By Lemma \ref{2.3} or Lemma \ref{2.4},
there are two points $p_1,q_1\in (p,q)$
satisfying the following properties:
\item{(1)} $q_1\in F({r})$ and $p_1\in F_0({r})$.
\item{(2)} $(p_1,q_1)\cap F_0(r)=\emptyset$
if $q_1>p_1$ and $(q_1,p_1)\cap F_0(r)=\emptyset$ if $q_1<p_1$.
\item{(3)} The length of the interval between the two points
$p_1$ and $q_1$ is less than or equal to $\delta/2$. \par
Since $p_1\in F_0({r})$,
we have $F(r)=F(p_1)$ by (L6). So we have
\item{($1'$)} $q_1\in F({p_1})$.
\item{($2'$)} $(p_1,q_1)\cap F_0(p_1)=\emptyset$
if $q_1>p_1$ and $(q_1,p_1)\cap F_0(p_1)=\emptyset$ if $q_1<p_1$.

We can repeat this argument replacing  $\{p,q\}$ by $\{p_1,q_1\}$.
Applying Lemma \ref{2.3} and Lemma \ref{2.4} inductively, we find
sequences $(p_n)$ and $(q_n)$ satisfying the following properties:
\item{(a)} $p_n$ lies in the interval beteen $p_{n-1}$ and $q_{n-1}$,  
and $q_n\in F({p_n})$.
\item{(b)} $(p_n,q_n)\cap F_0(p_n)=\emptyset$
if $q_n>p_n$ and $(q_n,p_n)\cap F_0(p_n)=\emptyset$ if $q_n<p_n$.
\item{(c)} The length of the interval between the two points
$p_n$ and $q_n$ is less than or equal to $\delta/2^n$.

It follows from Lemma \ref{2.1} and Lemma \ref{2.2} that
$$
\pi(F(p_n))\subset \pi(p_{n-1},q_{n-1}).
$$
In particular, the length of $\pi(F(p_n))$ is less than  
$\delta/2^{n-1}$.
We set
$$
y=\lim p_n=\lim q_n.
$$
The limit $y$ lies between $p_n$ and $q_n$ for all $n$.

We will now prove that $\pi(F(y))=\{\pi(y)\}$. Suppose that $\pi(F(y))$
does not only consist of $\pi(y)$.
Then there is a point $z\in F(y)$ such that
$Ty>z>y$.  For sufficiently large $n$, we either have
$$
Tp_n>z>q_n>y>p_n
$$
or
$$
Ty >Tq_n>z>p_n>y.
$$
In both cases (L4) implies that $F(y)=F(p_n)$.
In particular $y\in F_0(p_n)$, which contradicts
$(q_n,p_n)\cap F_0(p_n)=\emptyset$.
Thus we can conclude that $\pi(F(y))=\{\pi(y)\}$, which implies that
$y$ is a positive clean inflection point.
This is a contradiction. Hence there is a positive clean inflection
point $s$ in $(p,q)$. By Lemma \ref{2.1},
we have $\pi(F(s))\subset \pi((p,q))$.
\end{proof}

By reversing the orientation of $S^1$,
Proposition \ref{2.5} implies the following

\begin{proposition}\label{2.6}
Let $p\in S^1$. Suppose that $q\in F(p)\cap (Tp,p)$
and $(q,p)\cap F_0(p)=\emptyset$.
Then there exists a positive clean inflection point $s$ in $(q,p)$
such that $\pi(F(s))\subset \pi((q,p))$.
\end{proposition}

\begin{corollary}\label{2.7} Let $p\in S^1$.
Suppose that $q\in F(p)\cap (p,Tp)$
and $q\not\in F_0(p)$. Then there exists a positive clean
inflection point $s$ in $(p,q)$ such that $\pi(F(s))\subset \pi((p,q))$
and $F(s) \cap F_0(p)=\emptyset$.
\end{corollary}

\begin{proof}
We set
$$
p'={\rm sup}F_0(p).
$$
Since $q\not\in F_0(p)$ and $F_0(p')=F_0(p)$, we have
$$
q>p'\ge p,\qquad (p',q)\cap F_0(p')=\emptyset.
$$
Applying Proposition \ref{2.5} to the pair $(p',q)$,
we find a positive clean inflection
point $s$ in $(p',q)\subset (p,q)$.
We have $F(s) \cap F_0(p)=\emptyset$
since $\pi(F(s))\subset \pi((p',q))$.

\end{proof}

Similarly we get the following corollary.

\begin{corollary} \label{2.8} Let $p\in S^1$.
Suppose that $q\in F(p)\cap (Tp,p)$ and $q\not\in F_0(x)$.
Then there exists a positive clean inflection point$s$ in $(q,p)$
such that $\pi(F(s))\subset \pi((q,p))$
and $F(s) \cap F_0(p)=\emptyset$.
\end{corollary}

Applying Corollary \ref{2.7} and Corollary \ref{2.8}, we get
the following:

\begin{corollary} \label{2.9}
Suppose that $q\in F(p)$ satisfies $q\ne Tp$
and $q\not\in F_0(p)$. Let $J$ be the open
interval bounded by $p$ and $q$.
Then there exists a positive clean inflection point $s$
in $J$ such that $\pi(F(s))\subset \pi(J)$
and $F(s) \cap F_0(p)=\emptyset$.
\end{corollary}

Theorem B in the introduction is a consequence of the following theorem
if the curve $\gamma$ meets a line in at most finitely many components.

\begin{theorem} \label{2.10}
Let $\{F(p)\}_{p\in S^1}$ be an intrinsic
line system.
Then there exist three positive clean inflection points
$s_1,s_2,s_3$ in $S^1$ such that $s_2\in (s_1,Ts_1)$
and  $s_3\in (Ts_1,s_1)$.
Moreover, the sets $F(s_1),F(s_2),F(s_3)$
are mutually disjoint.
\end{theorem}

\begin{proof}
Take a point $p$ which is not a clean inflection point.
Then there exists a point $q\in  F(p)$
such that $q\not\in F_0(p)$.  By Corollary \ref{2.9},
there is a clean inflection point  $s_1$ between $p$ and $q$.
By (L5), we have $\pi(F(Ts_1))\ne \pi(F_0(Ts_1))$.
Then there exists
a point $u\in (s_1,Ts_1)$ such that $u\in F(Ts_1)$
but $u\not\in F_0(Ts_1)$.
Then by Corollary \ref{2.9}, we find a clean inflection
point $s_2$ on $(u,Ts_1)\subset (s_1,Ts_1)$.
Notice that $Tu\in F(Ts_1)$ and  $Tu\not\in F_0(Ts_1)$. Hence  we find  
another
positive clean inflection
point $s_3$ on $(Ts_1,Tu)\subset (Ts_1,s_1)$ by Corollary \ref{2.9}.
The sets $F(s_3)$ and $F(s_2)$ are disjoint
since $F(s_2)\subset (u,Ts_1)$ and $F(s_3)\subset (Ts_1,Tu)$.

Suppose that $F(s_2)\cap F(s_1)\ne \emptyset$.
Since $F(s_2)=F_0(s_2)$ and $F(s_1)=F_0(s_1)$, we have
$F(s_2)=F(s_1)$ by (L6). Then $Ts_1\in F(s_2)$
contradicting $F(s_2)\subset (u,Ts_1)$.
  Thus $F(s_2)\cap F(s_1)= \emptyset$.
Similarly we show $F(s_3)\cap F(s_1)= \emptyset$.
\end{proof}

Until now, we have assumed that
  $\gamma$ meets a line in
at most finitely many components.
We now prove Theorem B in the general case using
that such curves are generic in the set of anti-convex curves. In the  
proof
we will need that the curve $\gamma$ is $C^2$. So far we only used that
it is $C^1$.

\medskip
\noindent
{\it Proof of Theorem B}.
Let $\gamma$ be an arbitrary anti-convex curve on $P^2$ that
we assume to be $\pi$-periodic, that is
$\gamma(t)=\gamma(t+\pi)$ for $t\in \R$.
A point $p\in \R^3\setminus \{0\}$ determines a point $[p]$
in $P^2$, where $[p]$ denotes the line in $\R^3$ spanned by $p$.
There is an $\pi$-antiperiodic $C^2$-regular map $F:\R\to \R^3$
such that $$
\gamma(t)=[F(t)]\in P^2
$$
where a map $F(t)$ is called {\it $\pi$-antiperiodic} if
it satisfies $F(t+\pi)=-F(t)$ for all $t\in \R$.
The map $F$ has the Fourier series expansion
$$
F(t)=a_0+\sum_{n=1}^\infty\biggl(
a_n \cos (2n+1)t + b_n \sin (2n+1)t
\biggr ),
$$
where $a_0,a_1,b_1,...$ are all vectors in $\R^3$
and this series converges uniformly to $F(t)$.
We set
$$
F_N(t)=a_0+\sum_{n=1}^N\biggl(
a_n \cos(2n+1)t + b_n \sin (2n+1)t
\biggr ).
$$
One can easily show that
$\gamma_N(t)=[F_N(t)]$ is also anti-convex regular curve for  
sufficiently
large $N$ since $\gamma$ is $C^2$. We set
$$
\hat \gamma_N(t)=\frac{F_N(t)}{|F_N(t)|}:\R\to S^2.
$$

By Theorem \ref{2.10}, there exists
three positive clean  inflection points $s_1(N),s_2(N),s_3(N)$
on $\hat \gamma_N(t)$
such that
$$
0\le s_1(N)<s_2(N)-\pi<s_3(N)<s_1(N)+\pi<s_2(N)<s_3(N)+\pi<2\pi.
$$
By taking a subsequence, we may assume that
$s_j(N)$ converges to $s_j$ for $j=1,2,3$.
Since $\hat \gamma$ is not a great circle,
clean positive inflection points do not accumulate to
clean negative inflection points.
Thus we have
$$
0\le s_1<s_2-\pi<s_3<s_1+\pi<s_2<s_3+\pi<2\pi.
$$
These six points may not be clean flexes.
However, the tangent great circles at these six points
topologically cross $\hat \gamma$ exactly twice. Hence the
corresponding tangent lines of $\gamma$ only cross $\gamma$ once.
\qed

\section{Further properties of intrinsic line  
systems}\label{sec:3}\medskip

In this section we  derive some properties
of intrinsic line systems,
which will be used in the next section to
prove Theorem A in the introduction.
Throughout this section we will assume that an
intrinsic line system $\{F_p\}_{p\in S^1}$ is given.

For a point $p\in S^1$, we set
\begin{align*}
Y(p)&:=F(p)\setminus (F_0(p)\cup TF_0(p)), \\
Y^+(p)&:=Y(p)\cap [p,Tp], \qquad
Y^-(p):=Y(p)\cap [Tp,p], \\
F^+(p)&:=Y^+(p)\cup F_0(p), \qquad
F^-(p):=Y^-(p)\cup T(F_0(p)).
\end{align*}

\begin{figure}[htb]
   \begin{center}
         \includegraphics[width=4.2cm]{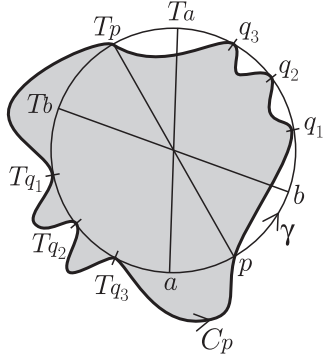}
\caption{Definition of $Y(p)$.}
\label{Yp}
   \end{center}
\end{figure}

For example, in the case of Figure \ref{Yp},
we have
$$
F_0(p)=\{p\},\qquad Y^+(p)=\{q_1,q_2,q_3\},
\qquad Y(p)=\{q_1,q_2,q_3,Tq_1,Tq_2,Tq_3\}.
$$

\begin{definition}\label{admissible}
An open interval $(a,b)$ is said to be {\it admissible} if
$b\in (a,Ta)$ and there are no positive
clean inflection points in $(a,b)$.
\end{definition}

\begin{figure}[htb] \label{defmuPm}
  \begin{center}
         \includegraphics[width=4cm]{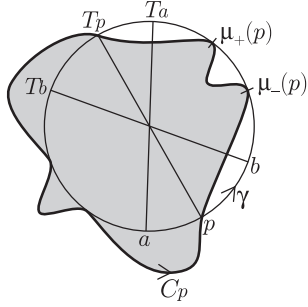}
\caption{Definition of $\mu_{\pm}(p)$.}
\label{cc4}
  \end{center}
\end{figure}

Let $(a,b)$ be an admissible interval.
Then $Y^+(p)$ is non-empty for all $p\in (a,b)$.
So we set (See Figure 13)
$$
\mu_-(p):=\inf_{(p,Tp)}Y^+(p),\quad
\mu_+(p):=\sup_{(p,Tp)}Y^+(p)
$$
for $p\in (a,b)$. 
For example, 
$$
\mu^{-}(p)=q_1,\qquad
\mu^{+}(p)=q_3
$$
holds in the case of Figure \ref{Yp}.
Moreover, we set
\begin{align*}
\mu_-(a):
&=
\begin{cases}
\inf_{[a,Ta]}{Y^+(a)}  & \mbox{if $a$ is not a positive
clean inflection point}, \\
\inf_{[a,Ta]}{TF_0(a)}  & \mbox{if $a$ is a positive
clean inflection point},
\end{cases} \\
\mu_+(b):
&=
\begin{cases}
\sup_{[b,Tb]}{Y^+(b)} & \mbox{if $b$ is not a positive
clean inflection point}, \\
\sup_{[b,Tb]}{F_0(b)} & \mbox{if $b$ is a positive
clean inflection point}.
\end{cases} \\
\end{align*}
Figure 14 explains  the definitions 
of $\mu_-(a)$ and $\mu_+(b)$ when
$a$ and $b$ are clean inflection points and
neither $F_0(a)$ nor $F_0(b)$ reduces a point.

\begin{figure}[htb]
   \begin{center}
         \includegraphics[width=4.2cm]{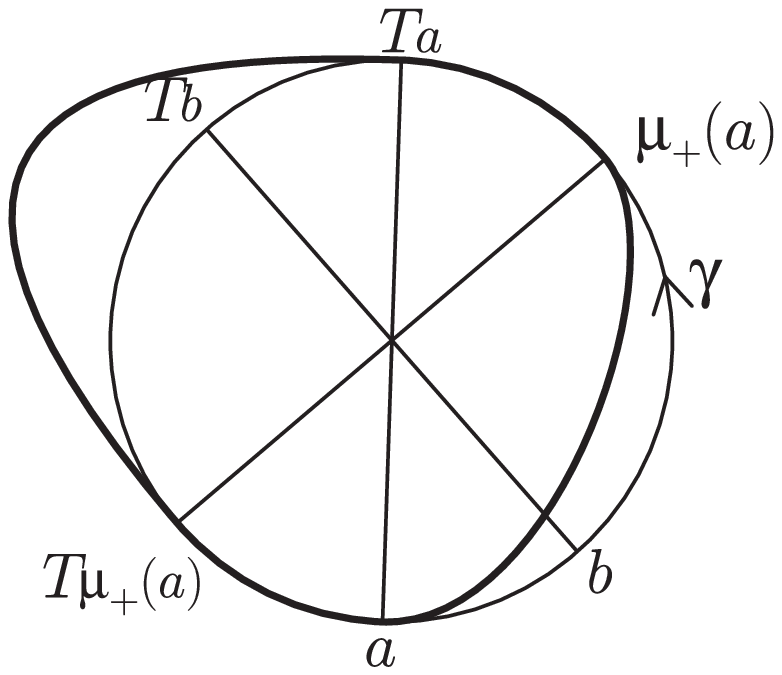}
\qquad
\includegraphics[width=3.8cm]{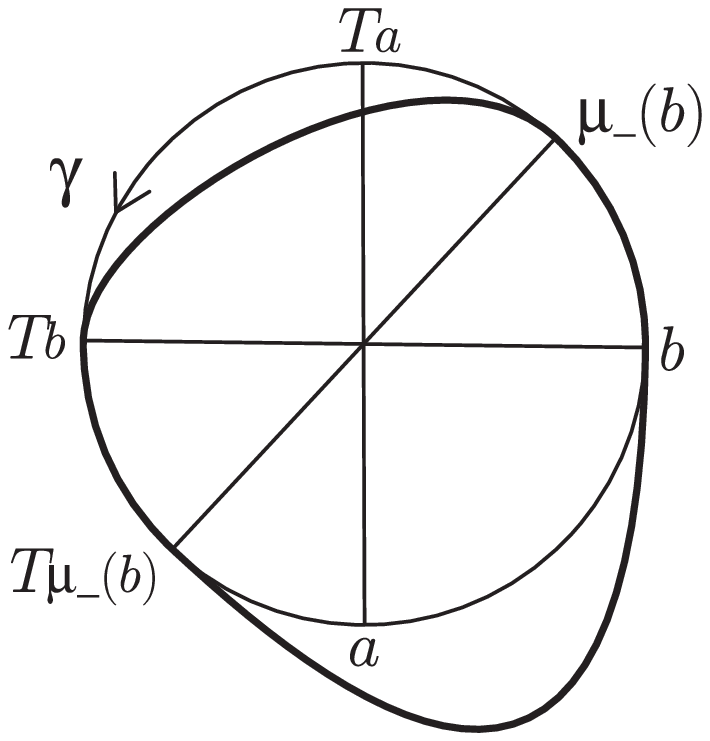}
\caption{Definitions of $\mu_-(a)$ and $\mu_+(b)$.}
\label{figs:edge}
   \end{center}
\end{figure}

These definitions have analogues in  the theory of
intrinsic circle system; see p.~190 in \cite{Umehara} by Umehara.
The results  in this section correspond to Lemma 1.3, 
Theorem 1.4 and Theorem 1.6 in \cite{Umehara}.
The left and the right of Figure \label{figs:edge}
correspond to the definition of $\mu_-(a)$ and
$\mu_+(b)$ when $a,b$ is positive 
clean inflection points, respectively.

\begin{remark}\label{rem:duality}
Let $S^1_{rev}$ be the $1$-dimesional sphere 
$S^1$ with the reversed orientation.
Then $\{F_p\}_{p\in S^1_{rev}}$ gives another 
intrinsic line system. 
An admissible interval $(a,b)$ of
$\{F_p\}_{p\in S^1}$
corresponds to the admissible interval
$(b,a)$ of $\{F_p\}_{p\in S^1_{rev}}$, 
and $\mu_-(p)$ ($p\in (a,b)$) with respect to 
$\{F_p\}_{p\in S^1}$ coincides
with $\mu_+(p)$ with respect 
to $\{F_p\}_{p\in S^1_{rev}}$. 
\end{remark}

\begin{lemma}\label{range0}
Let $(a,b)$ be an admissible interval.
Then  we have the inequalities
$$
b\le  \mu_+(p) < Ta
$$
for all $p\in (a,b]$
and
$$
b<  \mu_-(p) \le  Ta
$$
for all $p\in [a,b).$
\end{lemma}

\begin{proof}
We first assume that $p\in (a,b)$.
Then  $p$ is not a positive clean inflection point
and $Y^+(p)$ is non empty.
We fix $q\in Y^+(p)$  arbitrarily.
Then by Corollary \ref{2.9}, there is a positive
clean inflection point $r$ on $(p,q)$.
Since $(a,b)$ is an admissible arc, we have
$
q > r > b.
$
Suppose that
$$
(Tb>Tp) > q\ge Ta.
$$
Then we have
$$
b>p > Tq\ge a.
$$
Since $Tq\in Y^-(p)$, there is a positive clean inflection
point on $(Tq,p)\subset (a,b)$ by Corollary \ref{2.9}, which contradicts
the fact that $(a,b)$ is an admissible arc.
Thus we have $Ta> q$, which implies
$
q\in (b,Ta).
$
Since $q$ is arbitrary, we have
$$
b<\mu_-(p)\le \mu_+(p)<Ta
$$
for all $p\in (a,b)$.

\begin{figure}[htb]
  \begin{center}
         \includegraphics[width=4cm]{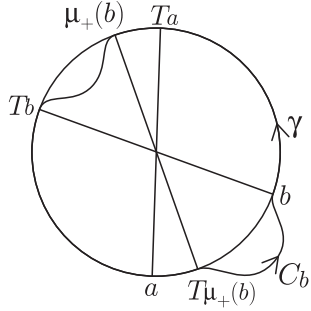}
\caption{The case $\mu_+(b)\ge Ta$.}
\label{cc5}
  \end{center}
\end{figure}

Next, we consider the case $q=b$.
If $b$ is not a positive clean inflection point,
$\mu_+(b)\in Y^+(b)$ and
the above arguments yield
$
b<\mu_+(b)<Ta.
$
So we assume $b$ is a positive clean inflection
point. Then
$
b\le \mu_+(b)
$
holds by definition.
Suppose now that
$
\mu_+(b)\ge Ta.
$ (See Figure \ref{cc5}.)
Then $T(\mu_+(b))\not\in F_0(b)$ and $\mu_+(b)\ne Tb$.
There is therefore a positive clean inflection point
between $(T(\mu_+(b)),b)$ by Corollary \ref{2.9},
which is a contradiction since $T(\mu_+(b))\in (a,b)$
and $(a,b)$ is admissible.
Thus we have
$
\mu_+(b)< Ta.
$

Finally, we consider the case $q=a$.
If $a$ is not a positive clean inflection point,
$\mu_-(a)\in Y^+(a)$ and
the above arguments yield $b<\mu_-(a)<Ta.$
So we assume $a$ is a positive clean inflection
point. Then $\mu_-(a)\le Ta$
holds by definition.
Suppose now that $\mu_-(a)\le b.$
(See Figure \ref{mu3}.)

\begin{figure}[htb]
  \begin{center}
         \includegraphics[width=4cm]{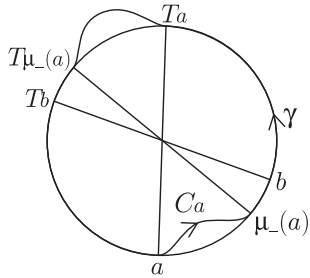}
\caption{The case $\mu_-(a)\le b$.}
\label{mu3}
  \end{center}
\end{figure}

Then $\mu_-(a)\ne Ta$.
Since $\mu_-(a)\not\in F_0(a)$,
there is a positive clean inflection point
between $(a,\mu_-(a))$ by Corollary \ref{2.9},
which is a contradiction since $T(\mu_-(a))\in (a,b)$
and $(a,b)$ is admissible.
Thus we have $b<\mu_+(b).$
\end{proof}

\begin{proposition}\label{range}
Let $(a,b)$ be an admissible interval. Then  we have the inequalities
\begin{align*}
(b\le )\mu_+(b)\le \mu_+(p), \\
\mu_-(p) \le \mu_-(a)(\le Ta)
\end{align*}
for all $p\in(a,b)$.
\end{proposition}

\begin{proof}
In the previous lemma, we already proved that
$$
b<  \mu_+(p)
$$
for all $p\in (a,b)$. Suppose now that
$
\mu_+(p)\in (b,\mu_+(b)).
$
Applying Lemma \ref{range0} to $(p,b)$,
we get
$
b\le \mu_+(b)<Tp.
$
Thus
$$
p<b<\mu_+(p)< \mu_+(b)(<Tp)
$$
holds.
Since $p,\mu_+(p)\in F^+(p)$,  we have
$F(b)=F(p)$ by (L4).
Thus  $b$ is like $p$ not a positive clean inflection and
$$
\mu_+(b)=\mu_+(p),
$$
contradicting the
the assumption $\mu_+(p)< \mu_+(b)$.
So we have
$
\mu_+(p)\ge  \mu_+(b).
$

By Lemma \ref{range0}, we have
$
\mu_-(p)<Ta.
$
Now we suppose
$$
\mu_-(a)< \mu_-(p)< Ta.
$$
Applying Lemma \ref{range0} to $(a,p)$,
we get
$
p<\mu_-(a).
$
Thus
$$
p< \mu_-(a) < \mu_-(p) < Ta (<Tp)
$$
holds.
Since $p,\mu_-(p)\in F^+(p)$, we have
$F(a)=F(p)$ by (L4).
Then  $a$ is like $p$ not a positive clean inflection
point.
Thus we have
$
\mu_-(a)=\mu_-(p),
$
contradicting
the assumption $\mu_-(a)<\mu_-(p)$.
So we have
$
\mu_-(p)\le  \mu_-(a).
$
\end{proof}

\begin{corollary}\label{monotonicity}(Monotonicity Lemma)
Let $(a,b)$ be an admissible arc and $p,q\in (a,b)$.
Suppose that $p<q$. Then we have
$$
\mu_-(p)\ge \mu_-(q),\qquad \mu_+(p)\ge \mu_+(q).
$$
Moreover
$
\mu_-(p)> \mu_+(q)
$
holds when $F(p)\ne F(q)$
and $\mu_-(a)>\mu_+(b) $ if there are points $p$ and $q$ in $(a,b)$
such that $F(p)\ne F(q)$.
\end{corollary}

\begin{proof} The first two inequalities follow directly from
Proposition~\ref{range}.

We now prove that
$\mu_-(p)> \mu_+(q)$
when $F(p)\ne F(q)$.
Assume that $F(p)\ne F(q)$
and $\mu_-(p)\le  \mu_+(q)$.
By Proposition \ref{range} we have
$$
(a<)p<q<\mu_-(p)\le  \mu_+(q)<Ta.,
$$
which implies by (L4) that $F(p)= F(q)$, which is a
contradiction. Hence
$\mu_-(p)> \mu_+(q)$.

Finally we prove the inequlity $\mu_-(a)>\mu_+(b)$
under the assumption that there are points $p,q\in (a,b)$ such that
$p<q$ and $F(p)\ne F(q)$.
 From Proposition \ref{range} and the inequality we have just proved
follows that
$$
\mu_{+}(b)\le \mu_{+}(q)< \mu_-(p) \le \mu_-(a)
$$
which proves the claim.
\end{proof}

\begin{proposition}\label{s-cont}(Semi-continuity)
Let $(a,b)$ be an admissible arc.
Then
$$
\lim_{x\to a+0}\mu_-(x)=\mu_-(a), \qquad
\lim_{x\to b-0}\mu_+(x)=\mu_+(b) .
$$
\end{proposition}

\begin{proof}
We shall prove the first formula. The second formula
can be proved similarly.
(See Remark \ref{rem:duality}.)
When there is a point $p\in(a,b)$ such that $p\in F_0(a)$, the  
assertion is obvious.
So we may assume that $(a,b)\cap F_0(a)=\emptyset$.
Let $(r_n)$ be a strictly decreasing sequence in $(a,b)$ converging to
$a$.  There are points $p_n$ and $q_n$ in the interval $(a,r_n)$ such  
that
$F(p_n)\ne F(q_n)$ since otherwise the closed set $F(q)$ would contain
the interval $[a,r_n]$ for all $q\in (a,r_n)$ and it would follow that  
$[a,r_n]\subset F_0(a)$.
Hence by Proposition \ref{range} and Corollary
\ref{monotonicity} we have that
$$
\mu_+(r_1)< \mu_-(r_n)<\mu_-(r_{n+1})< \mu_-(a)
$$
holds. So the sequence $\mu_-(r_n)$ has a limit  $s$.
Since $\mu_-(r_n)\in F(r_n)$,  (L7) implies that
$$
s \in F(a).
$$
Since
$\mu_+(p_1)\le  q_n\le \mu_-(a)$,
we have
$
\mu_+(p_1)\le s\le \mu_-(a).
$
Since $(a,b)\cap F_0(a)=\emptyset$, we have that
$(a,\mu_-(a))$ is disjoint from
the set $F(a)$.
Thus we have
$
s= \mu_-(a)
$
since $s\in F(a)$.
\end{proof}

\begin{theorem}\label{intermediate}
Let $(a,b)$ be an admissible arc.
Then for any $q\in (\mu_+(b),\mu_-(a))$,
there exists a point $p\in (a,b)$ such that
$$
\mu_-(p)\le q \le \mu_+(p).
$$
\end{theorem}

\begin{proof}
We set
$$
B_q:=\{x \in (a,b)\,;\,
  \mu_+(x) \le q\}.
$$
By
Proposition \ref{s-cont} we have that
$\displaystyle\lim_{x\to b-0}\mu_+(x)=\mu_+(b)+0$.
Thus a point $x\in (a,b)$ sufficiently close to
$b$ belongs to $B_q$. Since $B_q$ is non-empty,
we can set
$$
p:=\inf_{[a,b]}(B_q).
$$
Since $\mu_-(a)>q$, we have $p\in (a,b)$.
By the definition of $p$, there exists a sequence $(r_n)$ in  $B_q$
such that $\displaystyle\lim_{n\to \infty}r_n=p+0$.
By definition of $B_q$, we have
$$
\mu_-(r_n)\le \mu_+(r_n) \le q.
$$
Since $\displaystyle\lim_{n\to \infty} \mu_-(r_n)=\mu_-(p)$
by Proposition \ref{s-cont}, we have
$$
\mu_-(p)\le q.
$$
On the other hand, let $(s_n)$ be a sequence such that
$\displaystyle\lim_{n\to \infty}s_n=p-0$.
By definition of $B_q$, we have
$q<\mu_+(s_n)$.
Since $\displaystyle\lim_{n\to \infty} \mu_+(s_n)=\mu_+(p)$,
we have $q\le \mu_+(p)$.
\end{proof}

\section{Double tangents}\label{sec4}\medskip

  We we will assume throughout this section that
$\gamma:P^1\to P^2$ is an anti-convex $C^1$-regular curve  whose
number $i(\gamma)$ of true inflection points
is finite. It follows from the last assumption that
{\it a line in $P^2$ meets the curve $\gamma$ in at most finitely
many components.}

\begin{lemma}\label{lem4-1}
Let $\gamma:P^1\to P^2$ be an anti-convex curve.
Suppose that $\gamma$ meets a line  $L$ in $\gamma(a)$ and $\gamma(b)$
and denote one of the closed intervals on $P^1$ bounded by $a$ and $b$  
by $[a,b]$.
Then one of the two closed line segments $L_1$ and $L_2$ on $L$ bounded
by  $\gamma(a)$ and $\gamma(b)$,
say $L_1$, has the property that
$\gamma([a,b])\cup L_1$ lies in an affine plane and $\gamma([a,b])\cup  
L_2$
is not homotopic to a point.
The curve $\gamma([a,b])\cup L_1$ bounds a contractible domain
having acute interior angles at $\gamma(a)$ and $\gamma(b)$ if it is  
free
of self-intersections.\end{lemma}

We call $L_1$ the {\it chord  with respect to the interval
$[a,b]$} and denote it by $\overline{\gamma(a)\gamma(b)}$.

\begin{proof}
We choose a point $c\not\in [a,b]$. Then there is a
line $L_c$ which meets $\gamma$ only in $\gamma(c)$.
Then $L_c$ meets $L$ in one point which we assume to be
on the  line segments  on $L$ bounded
by  $\gamma(a)$ and $\gamma(b)$ that we denote by $L_2$.
Then $\gamma([a,b])\cup L_1$ lies in an affine plane.

  Since $L$ is not null-homotopic,
either
$\gamma([a,b])\cup L_1$ or  $\gamma([a,b])\cup L_2$
is not null-homotopic. So $\gamma([a,b])\cup L_2$
is not homotopic to a point.

Assume $\gamma([a,b])\cup L_1$  is free
of self-intersection and
let $D$ denote the contractible domain in the affine plane bounded by  
$\gamma([a,b])\cup L_1$.
If its interior angle at $\gamma(a)$ or $\gamma(b)$ is
not acute, any line passing through the point meets $\gamma$,
which contradics the anti-convexity of $\gamma$.
\end{proof}

The following assertion is one of the fundamental properties of
anti-convex curves.

\begin{proposition}\label{prop4-2}
Let $\gamma:P^1\to P^2$ be an anti-convex curve.
Let $[a,b]$ be a closed interval on $P^1$
and
suppose $\gamma([a,b])$ meets a line $L$ in $A^2$
at
$$
a=t_1<t_2<\cdots <t_n=b.
$$
Then
$$
\gamma(t_1),\,\,\gamma(t_2),\dots,\,\, \gamma(t_n)
$$
lie on $\overline{\gamma(a)\gamma(b)}$ in this order.
\end{proposition}

\begin{proof} Assume that the claim is not true. Then there is a  
smallest $i$
such that $\gamma(t_i)$ lies on  $\overline{\gamma(a)\gamma(t_{i-1})}$.
Then any line passing through $\gamma(t_i)$ must
meet $\gamma((t_1,t_i))$, which contradicts
the anti-convexity of $\gamma$.
\end{proof}

By Lemma \ref{lem4-1}, $\gamma([a,b])$ and
the chord $\overline{\gamma(a)\gamma(b)}$  lie in
an affine plane $A^2$. We
define a new curve $\gamma_1:P^1\to P^2$ by setting
$$
\gamma_1(t):=
\begin{cases}
\displaystyle
\frac{\gamma(b)(t-a)+\gamma(a)(t-b)}{b-a} \quad & \text{for  }\; t\in  
[a,b], \\
\gamma(t) \quad & \text{for }\; t\not\in (a,b)),
\end{cases}
$$
which is the curve one gets by replacing $\gamma([a,b])$ by
$\overline{\gamma(a)\gamma(b)}$. Notice that the  the vector operations  
in the definition of $\gamma_1$
depend on the affine plane $A^2$. 
We call $\gamma_1$ the
{\it  reduction of $\gamma$ with respect
to the  interval $[a,b]$.}\medskip

An interval $[a,b]$ on $P^1$ is called an {\it inflection
interval} if $a$ is a true inflection point
and $\gamma([a,b])$ is the connected component
of $\gamma\cap L_a$, where $L_a$ is the tangent line of $\gamma$ at $a$.

\begin{definition}\label{def:4-3}
Let $\gamma:P^1\to P^2$ be an anti-convex curve.
A nonempty proper open subinterval $(a,b)$ on $P^1$ is called
an {\it double tangent interval} if
\begin{enumerate}
\item the chord $\overline{\gamma(a)\gamma(b)}$
is tangent to $\gamma$ at $\gamma(a)$ and $\gamma(b)$.
\item there is a point in $\gamma([a,b])$ which is not
contained in $\overline{\gamma(a)\gamma(b)}$,
\item  $[a,b]$ is not an inflection interval of $\gamma_1$ where  
$\gamma_1$ is  the reduction of $\gamma$ with respect
to the interval $[a,b]$.\end{enumerate}
\end{definition}

By Lemma \ref{lem4-1},
the following assertion is obvious.

\begin{corollary}\label{lem4-2b}
If $(a,b)$ is a double tangent interval of an
anti-convex curve $\gamma$, then
the orientations of the  tangent lines of $\gamma$ at $\gamma(a)$ and  
$\gamma(b)$
induce the same direction on
$\overline{\gamma(a)\gamma(b)}$.
\end{corollary}

\begin{remark} \label{rem4-5}
If $(a,b)$ is a double tangent interval, then the same cannot be true  
for
$(b,a)=P^1\setminus [a,b]$.
In fact, the reduction $\gamma_2$ of $\gamma$ with respect
to the interval $[b,a]$ has $[b,a]$ as an inflection
interval  which violates property  (3) in Definition \ref{def:4-3}.
This phenomenon is explained in Figure \ref{gamma1} where
the two sketches indicate the same curve $\gamma$ in
different affine planes.

\begin{figure}[htbt]
   \begin{center}
          \includegraphics[width=3.5cm]{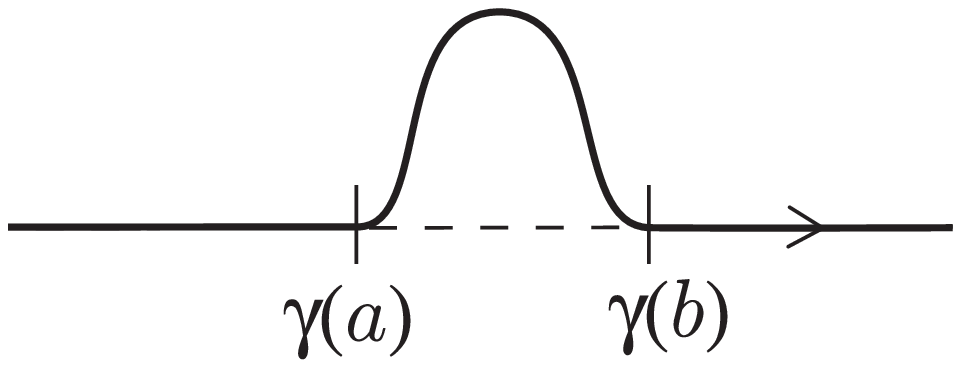}
\qquad
          \includegraphics[width=4cm]{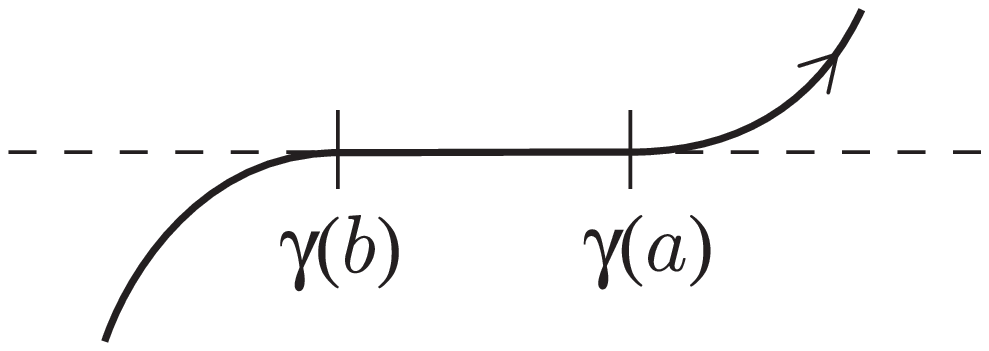}
\caption{ $\gamma$ in the different affine planes}
\label{gamma1}
   \end{center}
\end{figure}
\end{remark}

\begin{definition}\label{def4-5}
Let $\gamma:P^1\to P^2$ be an anti-convex curve.
Two double tangent intervals $(a_1,b_1)$
and $(a_2,b_2)$ are called {\it independent}
if they are disjoint or if the closure of one is contained in
the other.
\end{definition}

We now begin the proof of Theorem A in Introduction.

\begin{figure}[htbt]
  \begin{center}
         \includegraphics[width=4cm]{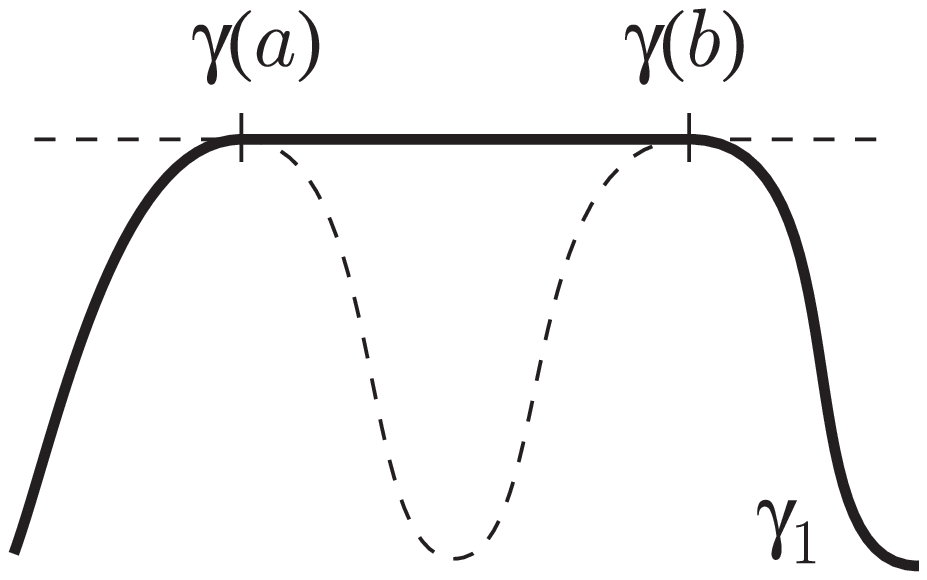}
\qquad\qquad
         \includegraphics[width=4cm]{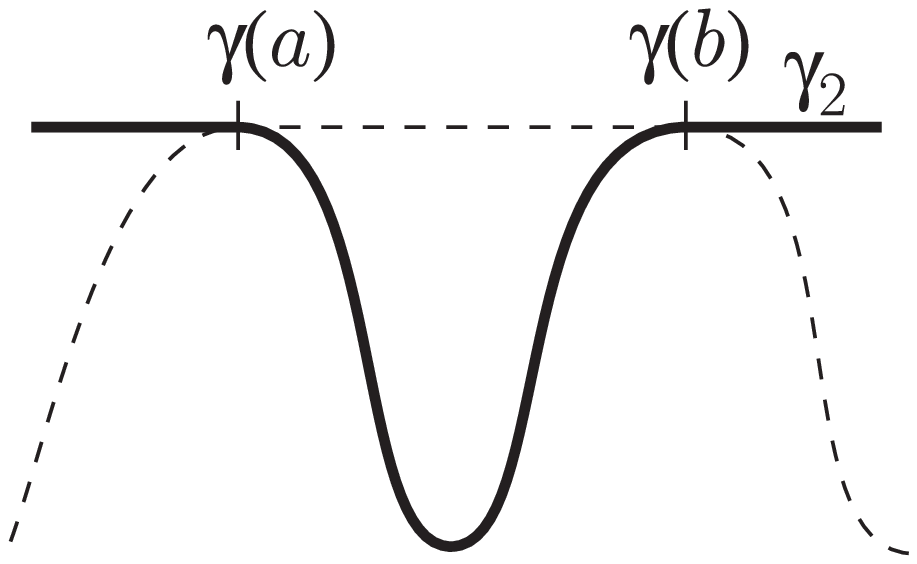}
\caption{$\gamma_1$ and $\gamma_2$}
\label{reduction}
  \end{center}
\end{figure}

\medskip
\noindent
{\it Proof of Theorem A.}
To prove   formula $(*)$, we will start with
a double tangent interval $(a,b)$ and introduce the
following reductions of $\gamma$.
We let $\gamma_1$ be
the reduction of $\gamma$ with respect
to the double tangent interval $[a,b]$
and we let
$\gamma_2$ be the reduction of
$\gamma$ with respect
to the  interval $[b,a]$; see Figure \ref{reduction}.

We now bring a couple of lemmas and propositions that will
be needed to finish the proof
of Theorem A.

\begin{lemma}\label{lem:simple}
The curves $\gamma_1$ and $\gamma_2$ are both
without self-intersections.
\end{lemma}

\begin{proof}
We will prove the claim  for $\gamma_1$.
Suppose $\gamma(P^1\setminus [a,b])$ meets the chord
$\overline{\gamma(a)\gamma(b)}$ at $\gamma(c)$.
By Proposition \ref{prop4-2} the points $\gamma(a)$,
$\gamma(b)$, $\gamma(c)$ must lie on the segment
$\overline{\gamma(a)\gamma(c)}$ in this order
since $a<b<c$. This is a contradiction. It follows that
$\gamma_1$ does not have self-intersections. One can
similarly prove that $\gamma_2$ does not have self-intersections.
\end{proof}

The following is a key to prove formula $(*)$.

\begin{proposition}\label{key1}
The curves $\gamma_1,\gamma_2$ are
both anti-convex and the identity
\begin{equation}\label{eq:inflection}
i(\gamma)= i(\gamma_1)+i(\gamma_2)-1
\end{equation}
holds.
\end{proposition}

\begin{proof}
We  first show that $\gamma_1$
is anti-convex. We may assume that $\gamma([a,b])$
lies in an affine plane $A^2$.
  For a point $x\in P^2$,
the pencil of lines passing through $x$ is a
a projective line in the dual space of $P^2$ that we denote by $P^1(x)$.
For a point $t\in P^1$,
we define a subset $\mathcal B_\gamma(t)$ of
$P^1(\gamma(t))$ such that
each line $L$ in $\mathcal B_\gamma(t)$ meets $\gamma$ only at $p$
and $L$ is transversal to the tangent line at $p$.
Since $\gamma(t)$ is an anti-convex curve,
$\mathcal B_\gamma(t)$ is non-empty for all $t\in P^1$.
One can easily prove that
  $\mathcal B_\gamma(t)$ is an open interval in $P^1(x)$.
We will  call $\mathcal B_\gamma(t)$ the {\it Barner set } of $\gamma$.

We have that $\mathcal B_\gamma(t)$ is contained in
the Barner set $\mathcal B_{\gamma_1}(t)$ of $\gamma_1$ for every  
$t\not\in [a,b]$,
since no line $L\in \mathcal B_\gamma(t)$ can meet the chord
$\overline{\gamma(a)\gamma(b)}$.
So it is sufficient to show that
$\mathcal B_{\gamma_1}(t)$ is not empty for $t\in (a,b)$.
Suppose $\gamma:P^1\to P^2$ meets the chord
$\overline{\gamma(a)\gamma(b)}$ at
$$
a=t_1<t_2<\cdots <t_n=b.
$$
By Proposition 4.2,
$$
\gamma(t_1),\,\,\gamma(t_2),\dots,\,\, \gamma(t_n)
$$
lie on $\overline{\gamma(a)\gamma(b)}$ in this order.

Suppose now that there exists a point $x\in  
\overline{\gamma(a)\gamma(b)}$
such that
the Barner set of $\gamma_1$ at $x$
is empty.
Then there exists a positive integer $i$ where $1\le i\le n-1$
such that $x\in \overline{\gamma(t_i)\gamma(t_{i+1})}$
and $x\ne \gamma(t_i),\gamma(t_{i+1})$.

\begin{figure}[htb]
  \begin{center}
         \includegraphics[width=6cm]{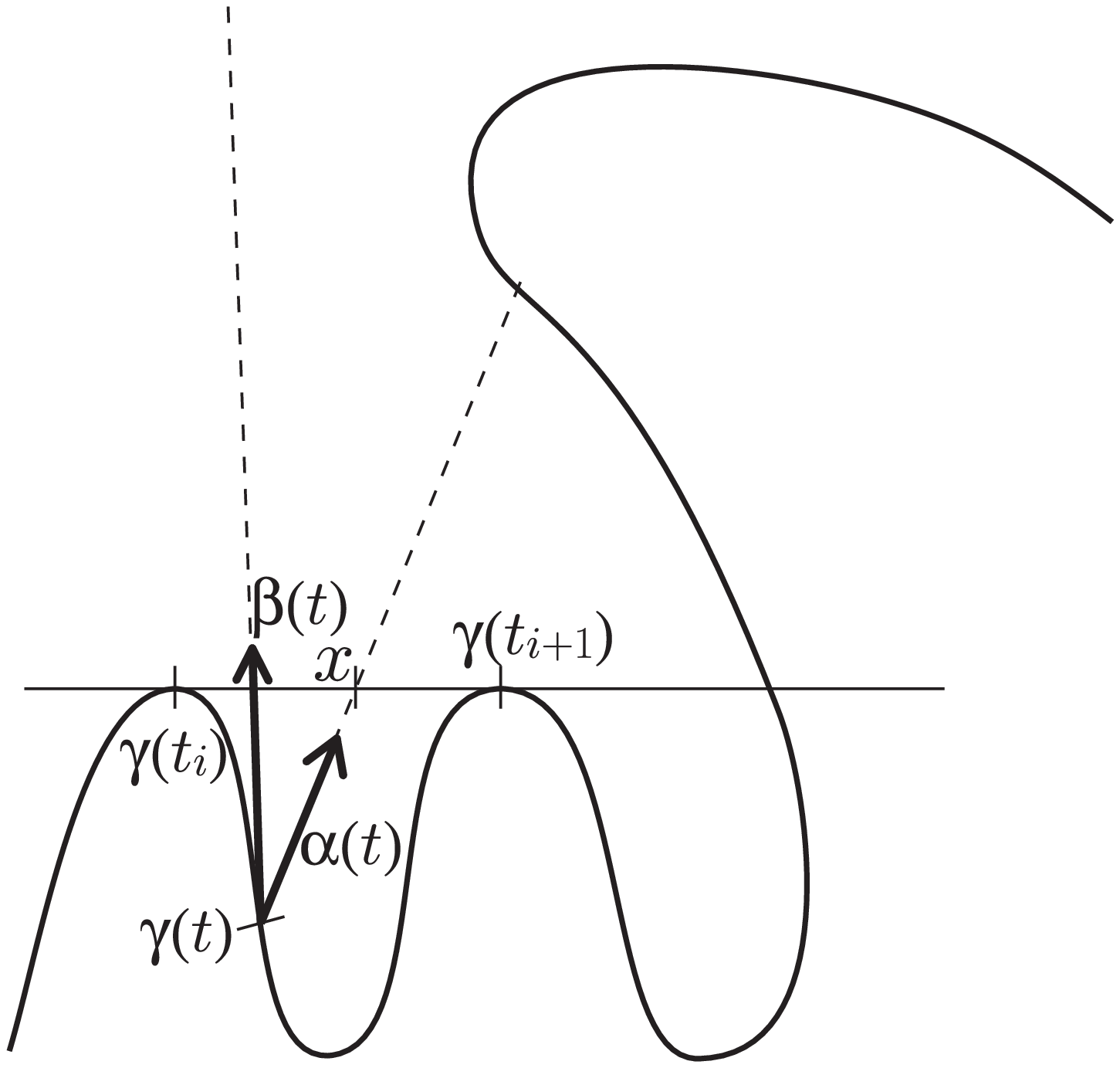}
\caption{}
\label{barner}
  \end{center}
\end{figure}

We now set
$$
I:=[t_i,t_{i+1}].
$$
In the following argument we work in $A^2$ that we equip with the  
orientation
such that
$\overline{\gamma(t_i)\gamma(t_{i+1})}$
lies on the left hand side of
$\gamma(I)$ as in Figure \ref{barner}.
We define continuous vector fields $\vect \alpha$ and $\vect \beta$
along  $\gamma|_I$  as follows:

\begin{enumerate}
\item Both $\vect \alpha(t)$ and $\vect \beta(t)$ point the left
hand side of $\gamma$ for every $t\in I$,
\item $\vect \alpha(t)$ lies on the chord $\overline{\gamma(t)x}$,
\item $\vect \beta(t)$ generates a line in $\mathcal B_{\gamma}(t)$.
\end{enumerate}
We set
\begin{align*}
I_L&:=\{t\in I\,;\, \mbox{$\vect \alpha(t),\vect \beta(t)$
is a positive frame}\}, \\
I_R&:=\{t\in I\,;\, \mbox{$\vect \alpha(t),\vect \beta(t)$
is a negative frame}\},
\end{align*}
that is $I_L$ (resp.~$I_R$) consists of those $t$  with the property  
that
the Barner direction $\vect\beta(t)$ is on the left of (resp.~right of)
of $\overline{\gamma(t)x}$.

Notice that $\vect \alpha(t)$ and $\vect \beta(t)$ are linearly  
independent
for all $t\in I$ since the Barner set of $\gamma_1$ at $x$ is empty.
Hence the sign of $\op{det}(\vect \alpha(t),\vect \beta(t))$ is either  
positive or
negative implying that either $I_L$ or $I_R$ is empty.
By Corollary \ref{lem4-2b} the tangent lines of $\gamma$ at $\gamma(a)$  
and
$\gamma(b)$ induce the same direction on  
$\overline{\gamma(a)\gamma(b)}$.
Hence it follows that
$
t_i \in I_L
$
and
$
t_{i+1} \in I_R
$,
and thus that neither $I_L$ nor $I_R$ is empty.
This is a contradiction and  we can conclude that
the Barner set of $\gamma_1$ at
a point $x\in \overline{\gamma(a)\gamma(b)}$
is not empty. This finishes the proof that $\gamma_1$ is anti-convex.
The proof that $\gamma_2$ is anti-convex is analogous.

Next we prove formula \eqref{eq:inflection}.
Let $I_1$ and $I_2$ be the number of independent inflection points of  
$\gamma$
on  $S^1\setminus [a,b]$ and $[a,b]$  respectively.
By definition, it is obvious that
\begin{equation}\label{eq:i2}
i(\gamma_2)=I_2+1.
\end{equation}
In fact $[b,a]$ is an additional inflection interval
on  $\gamma_2$. This phenomenon was explained in
Remark \ref{rem4-5} and Figure \ref{gamma1} above.
On the other hand,  we have
\begin{equation}\label{eq:i1}
i(\gamma_1)=I_1.
\end{equation}
By \eqref{eq:i2} and \eqref{eq:i1}, we hence have
$$
i(\gamma_1)+i(\gamma_2)=I_1+I_2+1=i(\gamma)+1,
$$
which proves  \eqref{eq:inflection}.
\end{proof}

\begin{corollary}\label{first}
  If $i(\gamma)=3$, then there are no
double tangent intervals on $\gamma$.
\end{corollary}

\begin{proof}
Suppose that there is a double tangent interval. Then
we can consider the  anti-convex curves $\gamma_1$ and $\gamma_2$
as in Proposition \ref{key1}.
Since both $i(\gamma_1)$ and $i(\gamma_2)$ are at least $3$ by Theorem  
\ref{2.10},
we have
$$
i(\gamma)= i(\gamma_1)+i(\gamma_2)-1
\ge 3+3-1=5
$$
which contradicts $i(\gamma)=3$.
\end{proof}

We are assuming in this section that the number  $i(\gamma)$ is finite.
This has a consequence for  number of elements in a set consisting of
independent  double tangent intervals
as the next corollary shows.

\begin{corollary}\label{infinity}
The number of elements in a set  of independent double tangent intervals
  is finite. \end{corollary}

\begin{proof} We assume that this number is infinite.
Let $n$ be an arbitrary positive integer.
Then we can find independent double tangent intervals
$
(a_1,b_1),  (a_2,b_2), \dots, (a_n,b_n).
$
We order the intervals such that $(a_i,b_i)$ does not contain  
$(a_j,b_j)$ for $i<j$.
We can associate to $(a_1,b_1)$ two
anti-convex curves $\gamma^{(1)}_1$ and  $\gamma^{(1)}_2$
as was done before Lemma~\ref{lem:simple}.
Then we use the same construction to associate to $(a_2,b_2)$ and  
$\gamma^{(1)}_1$ two
new anti-convex curves $\gamma^{(2)}_1$ and $\gamma^{(2)}_2$.
In this way we can get a finite sequence of pairs of anti-convex curves
$\gamma^{(k)}_1$ and $\gamma^{(k)}_2$ for $k=1,...,n$.
By Proposition \ref{key1} we have
$$
i(\gamma)=i(\gamma^{(n)}_1)-n+\sum_{k=1}^ni(\gamma^{(k)}_2)
$$
Since $i(\gamma^{(k)}_1), i(\gamma^{(k)}_2)\ge 3$, we
have $i(\gamma)\ge 3-n+3n=3+2n$.
Since $n$ is arbitrary, this contradicts the fact that $i(\gamma)$
is finite.\end{proof}

The proof of  the next proposition relies on the
results of Section \ref{sec:3}.

\begin{proposition}\label{key2}
If there are no double tangent intervals on $\gamma$, then
$i(\gamma)=3$ holds.
\end{proposition}

Let $\hat \gamma:S^1\to S^2$ be the lift of $\gamma$ to a closed curve  
on $S^2$.
We will need the following lemma in the proof of the poposition.

\begin{lemma}\label{interval}
Let $(a,b)$ be an admissible interval on $S^1$ in the sense of
Definition \ref{admissible}.
Suppose that there are no double tangent intervals on $\gamma$.
Then there are no true inflection points on $(\mu_+(b),\mu_-(a))$.
\end{lemma}

\begin{proof}
Let $\{F(p)\}_{p\in S^1}$ be the intrinsic line system
associated to the lift $\hat \gamma$.
Suppose that there is a true inflection point
$c\in (\mu_+(b),\mu_-(a))$.
By Theorem \ref{intermediate},
there exists a point $p\in (a,b)$, such that
$$
\mu_-(p)\le c \le \mu_+(p).
$$
Since $c$ is a true inflection point,
the limiting great circle $C_p$ cannot
pass through $\hat\gamma(c)$.
This implies that there is a double tangent interval on $\gamma$.
This contradiction proves the claim.
\end{proof}

\medskip
\noindent
{\it Proof of Proposition \ref{key2}}.
By Theorem \ref{2.10}, there are at least three
positive clean inflection intervals $[a_1,a_2]$, $[b_1,b_2]$ and
$[c_1,c_2]$ on $S^1$ some of which may of course reduce to points.
We assume that
$$
a_1\le a_2<b_1\le b_2<c_1\le c_2
$$
and that there are no positive clean inflection points
on $(a_2,b_1)$ and $(b_2,c_1)$.

By Lemma \ref{interval}, there are no inflection points on
$(c_2,Tb_1)$
since $(b_2,c_1)$ is an admissible arc
and $\mu_+(c_1)=c_2,\,\, \mu_-(b_2)=Tb_1$.
Since $\pi((c_2,Tb_1))=\pi((Tc_2,b_1))$,
there are no inflection points on
\begin{equation}\label{eq:C}
A:=(c_2,Tb_1)\cup (Tc_2,b_1).
\end{equation}

There are also no positive clean inflection points
on $[a_2,b_1]$.
Applying Lemma \ref{interval} to the interval
$(a_2,b_1)$, we conclude that
there are no inflection points on
\begin{equation}\label{eq:B}
C:=(b_2,Ta_1)\cup (Tb_2,a_1).
\end{equation}
In particular, there are no positive clean inflection points
on
$$(c_2,a_1)=(c_2,Tb_1)\cup (Tb_1,Tb_2)\cup (Tb_2,a_1).$$
Applying Lemma \ref{interval} to the interval
$(c_2,a_1)$, we conclude that
there are no inflection points on
\begin{equation}\label{eq:A}
B:=(a_2,Tc_1)\cup (Ta_2,c_1).
\end{equation}
Now it follows from  \eqref{eq:C}, \eqref{eq:B}, \eqref{eq:A} that
there are no inflection points on
$$
S^1\setminus
([a_1,a_2]\cup  [Tc_1,Tc_2] \cup
[b_1,b_2]\cup [Ta_1,Ta_2] \cup
[c_1,c_2]\cup [Tb_1,Tb_2])
=A\cup B\cup C,
$$
and hence that $i(\gamma)=3$.
\qed

\medskip
We can now finish the proof of Theorem A.
We will let $\delta(\gamma)$ denote the number of elements in a maximal  
set
of independent double tangent intervals. The number $\delta(\gamma)$
is finite by Corollary  \ref{infinity}. It will follow from the proof  
that $\delta(\gamma)$
does not depend on the maximal set that was used to define it.

We shall prove formula $(*)$ by induction over $i(\gamma)$.
When $i(\gamma)=3$, then
$(*)$ holds since  $\delta(\gamma)=0$
by Corollary \ref{first}.
So we assume $(*)$ holds when
$i(\gamma)\le n-1$ and $n\ge 4$
and prove it for $i(\gamma)=n$.
Since $i(\gamma)\ge 4$,
there exists at least one double tangent interval $I=(a,b)$
by Proposition \ref{key2}.
There exist non-negative integers $i$ and $j$ such that
\begin{enumerate}
\item
$
I,I_1,,\dots, I_i,J_1, ,\dots, J_j
$
is a maximal family of independent  double tangent intervals.
\item $I_1, \dots, I_i$ are subets of $I$,
\item $J_1,\dots, J_j$ lie on $\subset P^1\setminus (a,b)$.
\end{enumerate}
Then we get two anti-convex curves $\gamma_1,\gamma_2$
with respect to $I=[a,b]$.
By the  induction assumption
$\delta(\gamma_1)$ and $\delta(\gamma_2)$ do not depend on the choice of
the set of independent double tangent intervals.
Since $I_1,\dots, I_i$
and $J_1, \dots, J_j$
are maximal sets of
independent double tangent intervals
on $\gamma_1$ and $\gamma_2$ respectively,
we have
$$
i+j+1=\delta(\gamma_2)+\delta(\gamma_2)+1.
$$
By \eqref{eq:inflection},
we have
$$
i(\gamma)-2(i+j+1)=
\biggl( i(\gamma_1)-2\delta(\gamma_1)\biggr )
+ \biggl( i(\gamma_2)-2\delta(\gamma_2)\biggr )-3.
$$
By the induction assumption,
$$
i(\gamma_1)-2\delta(\gamma_1)=i(\gamma_2)-2\delta(\gamma_2)=3.
$$
Thus we have
$$
i(\gamma)-2(i+j+1)=3,
$$
which implies that the number $i+j+1$
of the independent double tangent intervals
is independent of the choice of
$I,I_1, \dots, I_i,J_1,\dots, J_j$.
Thus we have $\delta(\gamma)=i+j+1$.
This finishes the proof.
\qed

\section{Anti-periodic functions and curves of constant width}\medskip

Before giving a proof of Theorem C in the introduction,
we shall explain some properties of periodic and anti-periodic
functions, which we will need.
We denote by $C^r({\bf R})$, where $r=1,2,\dots, \infty$,
the vector space of $r$ times continuously differentiable real valued  
functions
on $\R$.
We define the following finite dimensional
linear subspaces of $C^r(\R)$
\begin{align*}\label{eq:A}
{\mathcal A}_{2n+1}&:=
\biggl\{a_0+\sum_{k=1}^{n}\biggl(a_k \cos kt +
b_k \sin kt\biggl)\,;\,
a_0,a_1,\dots, a_n, b_1,\dots, b_n\in \R
\biggr\}, \\
{\mathcal A}_{2n}&:=
\biggl\{\sum_{k=1}^{n}\biggl(a_k \cos (2k-1)t +b_k \sin (2k-1)t
\biggl)\,;\,
a_1,\dots, a_n, b_1,\dots, b_n\in {\bf R}
\biggr\},
\end{align*}
where $n$ is any natural number.
Let $f$ be a $C^r$-function and $m\le r$ some natural number.
For each point $p$ on $\R$, there exists
a unique function $\varphi_p$ in ${\mathcal A}_m$ such that
$$
f(p)=\varphi_p(p),\,\, f'(p)=\varphi'_p(p),
\,\, f''(p)=\varphi''_p(p), \,\, \dots
\,\, , f^{(m-1)}(p)=\varphi^{(m-1)}_p(p),
$$
namely, $\varphi_p$ is the best
approximation of $f$
at $p$ in ${\mathcal A}_m$.
We call $\varphi_p$ the {\it osculating function of order $m$}
or
{\it ${\mathcal A}_m$-osculating function
at $p$}.
In general the $m$-th derivative $\varphi_p^{(m)}(p)$
at $p$ is not equal to $f^{(m)}(p)$. If however
$\varphi_p^{(m)}(p)=f^{(m)}(p)$ holds for $p$, then $p$ is called a
{\it flex of $f$ of order $m$}.

Consider the following differential operators on ${\bf R}$
\begin{align*}
L_{2n+1}&:=D(D^2+1)(D^2+2^2)\cdots (D^2+n^2), \cr
L_{2n}&:=(D^2+1)(D^2+3^2)\cdots (D^2+(2n-1)^2),
\end{align*}
where $D=d/dt$.
Then ${\mathcal A}_m$ is the kernel of the operator
$L_m$.
The following proposition is proved in the appendix of \cite{TU3},  
p.~135.

\begin{proposition}
A point $p$ is a flex of $f$ of order $m$ if and only if
$(L_mf)(p)=0$.
\end{proposition}

\medskip
A function $f:\R\to \R$ is called
{\it $\pi$-antiperiodic} if it satisfies $f(t+\pi)=-f(t)$.
We now introduce the concept of clean flexes
for $2\pi$-periodic and $\pi$-antiperiodic functions.

\begin{definition}
Let $m$ be an integer that we first assume to be odd.
Let $f$ be a $2\pi$-periodic $C^r$-function where
$r\ge m-1$.
A point $p$ is called a {\it clean flex of
order $m$} if the set of zeros of the difference function $f-\varphi_p$
is connected in ${\bf R}/2\pi{\bf Z}$.

We next assume that $m$ is an even integer and $f$ a  $\pi$-antiperiodic
$C^r$-function where $r\ge m-1$.
Then a point $p$ is called a {\it clean flex of
order $m$} if the set of zeros of $f-\varphi_p$
is connected in ${\bf R}/\pi{\bf Z}$.
\end{definition}

\begin{remark} One should notice that $f$ does only have to be  
$C^{m-1}$ in the definition
of a {\it clean flex} of order $m$, but we needed $C^m$-regularity in  
the definition of
a {\it flex} of order $m$. If $f$ is $C^m$, then a clean flex of order  
$m$
is a flex of order $m$  in the sense of the former definition. It is  
crucial for many of our
arguments to allow low differentiability. Example for this are  
constructions like the
reductions of curves with respect to an interval in Section \ref{sec4}.
\end{remark}

In \cite{TU3} the authors proved the following:

\medskip
\noindent
{\it Let $m$ be a positive odd integer and let
$f$ be a $2\pi$-periodic $C^{m-1}$-function.
Then $f$ has at least $m+1$ clean flexes of order $m$ in a period.}

\medskip
In \cite{TU3} only the case  where $m$ is odd is dealt with.
One can expect that
{\it a  generic $\pi$-antiperiodic $C^{m-1}$-function
has at least $m+1$ clean flexes in a period}. An indication for this is
the fact that such a function of class $C^m$
{\it has at least $m+1$ possibly not clean flexes}
as can be easily proved;  see the appendix of [TU3].
In this section we give an affirmative answer
for the problem  if $m=2$ and leave the general case as an
open question. Our result is stated in the next theorem.

\begin{theorem} \label{thm:A2}
Let $f:{\bf R}\to {\bf R}$ be
a $\pi$-antiperiodic $C^1$-function  not belonging
to ${\mathcal A}_2$.
Suppose that  the zero set of $f-\psi$ is discrete for
every $\psi$ in $\mathcal A_2$.
Then $f$ has at least three clean flexes
$t_1<t_2<t_3$ of order 2,
where $t_3<t_1+\pi$, with the property that
$f-\phi_{t_1}$ and $f-\phi_{t_3}$ change sign
from negative to positive in $t_1$ and $t_3$ respectively,
and
$f-\phi_{t_2}$ changes sign
from positive to negative in $t_2$.
\end{theorem}

\begin{figure}[h]
  \begin{center}
         \includegraphics[width=5.5cm]{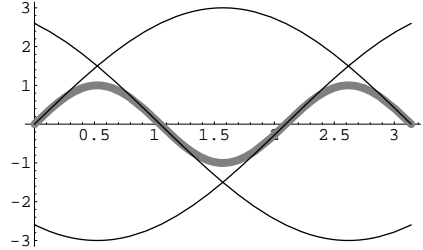}
  \end{center}
\caption{Three clean osculating functions for $\sin 3t$}
\label{sin3t}
\end{figure}

The theorem is optimal since $f(t)=\sin 3t$
has exactly three clean flexes at $t=0$, $\pi/3$, $2\pi/3$ in $[0,\pi)$;
see Figure \ref{sin3t}.
The theorem implies the well known existence of
three (usual) flexes of order $2$ which
can be proved by integration by parts; see \cite{GMO} and
the appendix of \cite{TU3}.

We start with some lemmas needed
to prove Theorem \ref{thm:A2}.
In the following $f$ will be a $\pi$-anticonvex $C^1$-function as
in the theorem.
For a point $p$ we define  a one-dimensional subspace $V_p$ of
$
{\mathcal A}_2=\{a \cos t + b\sin t\,;\, a,b\in \R\}
$
by setting
$$
V_p:=\{\psi\in {\mathcal A}_2\,;\, \psi(p)=f(p)\}.
$$
The osculating function $\varphi_p$ at $p$
belongs to $V_p$.
For a given $s\in {\bf R}$, there is a unique $\psi\in V_p$
such that $\psi'(p)=s$,
since ${\mathcal A}_2$ is the kernel of the operator $L_2$.
We will denote this function  by $\psi_{p,s}$.
Thus we may write $V_p=\{\psi_{p,s}\,;\,s\in {\bf R}\}$.
For sufficiently large $s$, the function $\psi_{p,s}$ has
the following properties:
\begin{enumerate}
\item[(1)] $\psi_{p,s}(t)$ is greater than $f(t)$ on $(p,p+\pi)$ and
\item[(2)] $\psi_{p,s}(t)$ is less than $f(t)$ on on $(p-\pi,p)$.
\end{enumerate}
Let $s_0$ be  the infimum over the set of real numbers $s$
such that $\psi_{p,s}$ satisfies (1) and (2)
and set
$$
\psi_p:=\psi_{p,s_0}.
$$
We will call $\psi_p$ the {\it limiting function of $f$ at $p$}.

\begin{lemma}\label{lem:maximal}
The limiting function $\psi_p(t)$ of $f$ at $p$
has the following three properties.
\begin{enumerate}
\item[{(a)}] $\psi_{p}(t)\ge f(t)$ for  $t\in (p,p+\pi)$,
\item[{(b)}] $\psi_{p}(t)\le f(t)$ for $t\in (p,p+\pi)$.
\item[{(c)}] If $\psi_p$ is not the ${\mathcal A}_2$-osculating function
of $f$ at $p$, then there exists a point $q$
on $(p,p+\pi)$ such that $\psi_p(q)=f(q)$.
\end{enumerate}
Conversely a function $\psi\in V_p$ satisfying $(a)$, $(b)$ and $(c)$
must coincide with the limitting function at $p$.
\end{lemma}

\medskip
\begin{proof}
The lemma is an analogue of Proposition \ref{1.2}
and Proposition \ref{1.3} and follows directly from the definition of  
$\psi_p$.
\end{proof}

Now we identify
$$
S^1={\R}/2\pi {\bf Z},\qquad
P^1={\R}/\pi {\bf Z},
$$
and denote by
$$
\pi:S^1\to P^1
$$
the canonical projection. We will consider $f$ and the limiting  
functions
$\psi_p$ as functions on $S^1$.
We now set
$$
F(p)=\{t\in S^1\,;\, f(t)=\psi_p(t)\}.
$$

\begin{lemma}\label{lem:clean}
A point $p$ is clean flex of order 2 if and
only if $F(p)$ consists of exactly two points.
\end{lemma}

\begin{proof}
If $p$ is a clean flex, then $\psi_p=\varphi_p$ by property (c).
Since the zero set of $f-\psi$ is discrete for
every $\psi$ in $\mathcal A_2$, $F(p)$ consists of exactly two points.
Conversely,  suppose $F(p)$ consists of exactly two points.
  Then
$\psi_p=\varphi_p$ by Lemma \ref{lem:maximal} and $p$ is a clean flex.
\end{proof}

\begin{proposition}\label{prop:axiom}
Let $f$ be as in Theorem \ref{thm:A2}.
Then the associated family of closed subset $\{F(p)\}_{p\in S^1}$
is an intrinsic line system.
\end{proposition}

\begin{proof} We have to show that properties (L1) to (L7) in  
Proposition \ref{1.4}
are satisfied.

(L1) is obvious. (L2) follows from the fact that
$f$ does not belongs to ${\mathcal A}_2$.
(L3) holds since $f$ is  $\pi$-antiperiodic.
(L4) follows from the fact that the functions in ${\mathcal A}_2$
have at most one zero on $[0,\pi)$.
(L6) is empty under our conditions. (L7) holds since the limit of a  
sequence
of limiting functions is a limiting function.
\end{proof}

\medskip
\noindent
{\it Proof of Theorem \ref{thm:A2}}.
We have associated an intrinsic line system $\{F(p)\}_{p\in S^1}$
to $f$ in Proposition \ref{prop:axiom}.
Now Theorem \ref{thm:A2}
implies that $f$ has at least three clean flexes of order 2 and it is  
easy to
see that they can be chosen as claimed in the theorem.
\qed

\begin{definition}\label{def:inflection}
A nonempty proper open subinterval $(a,b)$ on $P^1$ is called
an {\it $\mathcal A_2$-double tangent interval of $f$} if there is a
function $\varphi$ in $\mathcal A_2$ such that
\begin{enumerate}
\item the values of $f$ and  $\varphi$ coincide  in $a$ and $b$,
\item the derivatives of $f$ and $\varphi$ coincide in $a$ and $b$,
\item there is a point in $t\in (a,b)$ such that
$\varphi(t)\ne f(t)$,
\item the function $f-\varphi$ either has local maxima at both $a$ and  
$b$,
or it has local minima at both $a$ and $b$.
\end{enumerate}
\end{definition}

If $(a,b)$ is a double tangent interval, then the same cannot be true  
for
$(b,a)=P^1\setminus [a,b]$,
since condition (4) fails.
(If we consider $(a,b)$ to be
an interval of $\R$, then $P^1\setminus [a,b]$
corresponds to $(b,a+\pi)$.)
The function $\varphi$ in the definition of an
$\mathcal A_2$-double tangent interval
  is uniquely
determined.
We will call it the {\it double tangent function}
with respect to $(a,b)$.

\begin{definition}\label{def:independent}
Let $f:\R\to \R$ be an anti-periodic $C^1$-function.
Then two $\mathcal A_2$-double tangents $(a_1,b_1)$
and $(a_2,b_2)$ are said to be {\it independent}
if they are disjoint or if the closure of one is contained in
the other.
\end{definition}

Using the same method as in Section 4,
we get the following:

\begin{theorem}\label{A2main}
Let $f:{\bf R}\to {\bf R}$ be
a $\pi$-antiperiodic $C^1$-function  not belonging
to ${\mathcal A}_2$.
Suppose that  the zero set of $f-\psi$ is discrete for
every $\psi$ in $\mathcal A_2$.
Then the number $i(f)$ of flexes of order $2$, and
the number $\delta(\gamma)$ of elements in a maximal set of
independent ${\mathcal A}_2$-double tangent intervals
are both finite and $\delta(\gamma)$ is independent of
the choice of the maximal set of double tangent intervals.
Moreover
$$
i(\gamma)-2\delta(\gamma)=3,\qquad
$$
holds.
\end{theorem}

\begin{proof}
Let $(a,b)$ be an ${\mathcal A}_2$-double tangent interval
and $\varphi$ the corresponding double tangent function in ${\mathcal  
A}_2$.
Without loss of generality,
we may asuume that $0\le a< b<\pi$.
Then we set
$$
f_1(t):=
\begin{cases}
\varphi(t) \quad & \text{for  } t\in [a,b], \\
f(f)    \quad & \text{for  } t\in [0,\pi)\setminus [a,b],
\end{cases}
$$
and extend $f_1$ to $\R$ as a $\pi$-antiperiodic function.
Then $f_1$ is  a $C^1$-function that we call the {\it reduction of $f$
with respect to $[a,b]$.}
Similarly  we set
$$
f_2(t):=
\begin{cases}
f(t) \quad & \text{for  } t\in [a,b], \\
\varphi(t)    \quad & \text{for  } t\in [0,\pi)\setminus [a,b],
\end{cases}
$$
and extend $f_2$ to $\R$ as a $\pi$-antiperiodic function.
Then $f_2$ is a $C^1$-function that we call {\it the reduction of $f$
with respect to $[b,a]$.}  We now use the two functions $f_1$ and $f_2$
as the we used the reductions $\gamma_1,\gamma_2$ in
Section 4 to prove  Theorem A in the introduction by induction.
\end{proof}

Finally we come to the applications of Theorem \ref{thm:A2}
and Theorem \ref{A2main}
to convex curves of constant width in the Euclidean plane ${\bf R}^2$.

We first describe the connection between strictly convex curves and  
periodic
functions, where we understand under a  {\it strictly convex curve}  a  
convex curve with the
property that the tangent lines at different points are different.
Let  $o$ be a point in the open domain bounded
by a strictly convex $C^2$-regular curve $\gamma$ in $\bf R^2$.
For each $t\in [0,2\pi)$,
there is a unique tangent line $L(t)$ of the
curve which makes angle $t$ with the $x$-axis.
Let $h(t)$ be the distance
between
$o$ and the line $L(t)$. The $C^1$-function
$h$ is called the {\it supporting
function of
the curve $\gamma$ with respect to $o$.}
Set $\mathbf e(t)=(\cos t,\sin t)$ and
$\mathbf n(t)=(-\sin t,\cos t)$.
Then
$$
\gamma(t)=h'(t)\mathbf e(t)-h(t)\mathbf n(t)
$$
gives a parametrization of the curve $\gamma$.
The following lemma follows immediately.

\begin{lemma}\label{lem:circle}
Let $\gamma_1$ and $\gamma_2$ be two strictly convex
$C^2$-regular curves having a common point $o$ in their
interior, let $h_1(t)$ and $h_2(t)$ be
their supporting functions,
and let $\gamma_1(t)$ and $\gamma_2(t)$
be their parametrizations as above.
Then the difference $h_2(t)-h_1(t)$ does not depend on the
choice of the origin $o$.
In particular, if $\gamma_2$ is a circle,
then the point  $\gamma_1(t)$ lies in the interior of $\gamma_2$
if and only if $h_2(t)-h_1(t)>0$ holds.
  \end{lemma}

A convex curve has constant width $d$ if and only if
$h(t)+h(t+\pi)=d$ holds.

Let  $\gamma$ be a $C^2$-regular strictly convex closed curve
of constant width $d>0$ and $h$ its supporting function which is of  
class $C^1$.
The function $f_\gamma$ defined by
$$
f_\gamma(t)=h(t)-\frac{d}2
$$
is $\pi$-antiperiodic since $\gamma$ is of constant width.
If $\gamma$ is a circle of diameter $d$,
the supporting function
$\psi$
can be written as
$$
h(t)=\frac{d}2+b\cos t + c\sin t
$$
where $(c,-b)$ is the center of the circle.

For a point $p$ on a curve $\gamma$ of constant width $d$, there exists  
a
unique circle $\Gamma_p$
of width $d$ such that  $\Gamma_p$ is tangent to $\gamma$ at $p$,
that is $\Gamma_p$ and $\gamma$ meet at $p$ with multiplicity two.
Since $\Gamma_p$ is the best approximation of $\gamma$ at $p$ by a
circle of width $d$, we call $\Gamma_p$ the {\it osculating
$d$-circle at $p$.}
When $\Gamma_p$ meets $\gamma$ with multiplicity higher than two in $p$,
we call $p$ a {\it $d$-inflection point}.

\begin{proposition}\label{prop:osculating}
Let $\gamma$ be a $C^3$-regular convex
curve of of constant width $d$ and $h$
the supporting function of $\gamma$.
Then  the following three properties are equinvallent:
\begin{enumerate}
\item
a point $p=\gamma(t_0)$ is a $d$-inflection point,
\item
$h''(t_0)+h(t_0)={d}/{2}$,
\item the osculating  $d$-circle
$\Gamma_p$ at $p$ is an osculating circle in the usual sense,
that is, the curvature radius of $\gamma$ at $p$ is $d/2$.
\end{enumerate}
\end{proposition}

\begin{proof}
The supporting fuction $h$ is a $C^2$-function
because $\gamma$ is $C^3$-regular.
Since the radius of the osculating circle
of $\gamma$ at $t$ is given by
$r=h''(t)+h(t)$, the last two properties are
equivalent.
It is therefore sufficient to prove the equivalence of the first two
properties.

Let
$$
h_0=
\frac{d}2+b\cos t + c\sin t
$$
be the supporting function of a circle $\Gamma$.
Then $\Gamma$ is the $d$-osculating circle at $p$
if and only if
$$
h_0(t_0)=h(t_0),\quad h'_0(t_0)=h'(t_0).
$$
Moreover, $\Gamma$ and $\gamma$ meet
with multiplicity higher than two in $p$
if and only if the curvature radius of them coincide.
Since the radius of the osculating circle
of $\gamma$ at $t$ is given by
$r=h''(t)+h(t)$, the circle $\Gamma$ is a $d$-inflection point
if and only if
$$
r=h''(t_0)+h(t_0)=h''_0(t_0) + h_0(t_0)=\frac{d}2.
$$
This proves that the first two properties are equivalent.
\end{proof}

If $\Gamma_p\cap\gamma$ consists of exactly two connected
components, $p$ is a $d$-inflection point which we will call
a {\it clean $d$-inflection point}.

We can now prove Theorem C in the introduction as an application of  
\ref{thm:A2}.

\medskip
\noindent
{\it Proof of Theorem C.}
We first consider the special case that
  $\gamma$ meets  circles in at most finitely many points.
   As explained above,
the supporting function of $h$ can be written in the form
$$
h(t)=\frac{d}2+f_\gamma,
$$
where $f_\gamma$ is a  $\pi$-antiperiodic function.
Since we are assuming that $\gamma$ meets circles
in at most finitely many points, $f_\gamma-\psi$ has a discrete zero  
set for every
$\psi$ in $\mathcal A_2$.
Hence, by Theroem \ref{thm:A2}, the function $f_\gamma$ has at least  
three clean
positive flexes $t_1,t_2,t_3$ of order $2$ on the interval $[0,\pi]$.
By Proposition 5.1,
$f_\gamma''+f_\gamma$ vanishes at $t_1$, $t_2$, and $t_3$, that is
$$
h''(t)+h(t)=\frac{d}2
$$
holds for $t=t_1,t_2,t_3,t_1+\pi,t_2+\pi,t_3+\pi$.
By Proposition \ref{prop:osculating}, these six points
are $d$-inflection points.
Moreover, since these six points are clean flexes,
Lemma \ref{lem:circle} implies that
they clealy turn out to be six clean
$d$-inflection points. Notice that the corresponding
osculating $d$-circles meet $\gamma$ exactly twice in
$t_i$ and $t_i+\pi$. This finishes the proof of the special case.

Next we consider the general case in which
$\gamma$ can meet circles
infinitely many times.
We consider the Fourier series expansion
$$
h(t)=a_0+\sum_{n=1}^\infty
\biggl(
a_n \cos(2n+1)t + b_n \sin(2n+1)t
\biggr).
$$
of $h$ and  set
$$
h_N(t)=a_0+\sum_{n=1}^N
\biggl(
a_n \cos(2n+1)t + b_n \sin(2n+1)t
\biggr).
$$
Then the convex curve
$\gamma_N(t)$ with the supporting function $h_N(t)$
is a regular curve of constant width for  $N$ sufficiently large.
Now we apply  the same argumants as in the proof of Theroem B in  
Section 2
to find three distinct osculating $d$-circles as a limit of
those of $\gamma_N(t)$.
\qed
\medskip

One can of course  define the double $d$-tangent intervals
of a curve $\gamma$ of constant width as the double tangent
intervals of the corresponding function $f_\gamma$. We translate
this into geometric properties of $\gamma$ as follows.

A nonempty proper open interval $(a,b)$ of $S^1=\R/2\pi\mathbf Z$
is a {\it double tangent interval of $\gamma$} if there is a circle
$\Gamma$ which coincides with the
osculating $d$-circles at $\gamma(a)$ and $\gamma(b)$
and has the property that there is a $t$ in $(a,b)$ such that
$\gamma(t)\not\in\Gamma$. We assume furthermore that $\Gamma$ is locally
around $\gamma(a)$ and $\gamma(b)$ on the same side of $\gamma$.
Notice that $(a+\pi,b+\pi)$ is a double tangent interval if $(a,b)$ is  
such an interval.

Two double tangent intervals $(a_1,b_1)$ and $(a_2,b_2)$ are {\it  
independent} if
they are not antipodal on $S^1$ and if
they are disjoint or the closure of one is contained in the other.

Since  $\mathcal A_2$-double
tangent intervals correspond to the $d$-double tangent intervals,
Theorem \ref{A2main} implies the following theorem.

\begin{theorem}\label{A2cor}
Let $\gamma$ be a convex $C^3$-regular curve
of constant width $d$.
Suppose that the curve meets  circles
in at most finitely many points.
Then the number $i(\gamma)$ of independent
$d$-inflection points  and
the number $\delta(\gamma)$ of elements in a maximal set of
independent $d$-double tangent intervals
are both finite and $\delta(\gamma)$ is independent of
the choice of the maximal set used to define it.
Moreover we have the equation
$$
i(\gamma)-2\delta(\gamma)=3.
$$
\end{theorem}

%

\end{document}